\documentclass[10pt]{article}
\usepackage{amsmath, amssymb, amscd, eucal}
\usepackage{epsfig}
\usepackage{graphicx}
\newtheorem{theorem}{Theorem}
\newtheorem{corollary}{Corollary}
\newtheorem{lemma}{Lemma}

\begin{document}

\title{Computing and Sampling Restricted Vertex Degree Subgraphs and
Hamiltonian Cycles} \author{
  by \textsc{Scott Sheffield}
}
\date{}
\maketitle

\begin{abstract} Let $G=(V,E)$ be a bipartite graph embedded in a plane
(or $n$-holed torus).  Two subgraphs of $G$ differ by a {\it
$Z$-transformation} if their symmetric difference consists of the boundary
edges of a single face---and if each subgraph contains an alternating set
of the edges of that face. For a given $\phi: V \mapsto \mathbb Z^+$,
$S_{\phi}$ is the set of subgraphs of $G$ in which each $v\in V$ has
degree $\phi(v)$.  Two elements of $S_{\phi}$ are said to be adjacent if
they differ by a $Z$-transformation.  We determine the connected
components of $S_{\phi}$ and assign a {\it height function} to each of
its elements.

If $\phi$ is identically two, and $G$ is a grid graph, $S_{\phi}$ contains
the partitions of the vertices of $G$ into cycles.  We prove that we can
always apply a series of $Z$-transformations to decrease the total number
of cycles provided there is enough ``slack" in the corresponding height
function.  This allows us to determine in polynomial time the minimal
number of cycles into which $G$ can be partitioned provided $G$ has a
limited number of non-square faces.  In particular, we determine the
Hamiltonicity of polyomino graphs in $O(|V|^2)$ steps.  The algorithm
extends to $n$-holed-torus-embedded graphs that have grid-like properties.  
We also provide Markov chains for sampling and approximately counting the
Hamiltonian cycles of $G$.

\end{abstract}

\maketitle
\thispagestyle{empty}

\section{Introduction}

Let $G=(V,E)$ be a planar or $n$-holed-torus-embedded polygon graph (i.e.,
a graph whose faces are simple polygons) and let $\delta$ be a function on
the oriented edges of $G$.  Then the {\it height function space} $\mathcal
H = (G, \delta)$ is the set of (real or integer valued) functions $\phi$
on the face set $F$ of $G$ satisfying difference restrictions $\phi(f_2) -
\phi(f_1) \leq \delta(e)$ whenever $e$ is an edge shared by $f_1$ and
$f_2$ oriented with $f_1$ on the left.

If $G$ is planar, $\mathcal H$ is in one-to-one correspondence with the
set of closed electric flows $\omega$ on $G$ satisfying $-\delta(e^{-1}
\leq \omega(e) \leq \delta(e)$ for each oriented edge $e$.  This allows us
to convert questions about randomly chosen network flows or corresponding
objects (e.g., perfect matchings on bipartite graphs) into questions about
random height functions.

One advantage of this conversion is that recently introduced perfect
sampling techniques, including {\it Coupling from the Past} (see
\cite{PW}) and {\it Fill's algorithm} (see \cite{F}), apply especially
well to height function spaces.  Another is that large deviations
principles and asymptotic limiting shapes can sometimes be derived for
random height function spaces in a way that sheds light on local behavior.  
(See \cite{CEP} and \cite{CKP} for numerous references to examples
involving domino tilings.)

This paper uses height function theory to solve two other problems that
have been dealt with in the literature.  The first of these is {\it
$Z$-transformation connectedness}.  The problem here is as follows.  Let
$G$ be a bipartite, planar (or $n$-holed torus embeddable) graph
$G=(V,E)$, and $\phi$ a function $\phi:V \mapsto \mathbb Z^+$.  Let
$S_{\phi}$ be the set of all subgraphs $H$ of $G$ satisfying $deg_H(v) =
\phi(v)$ at every vertex $v$.

Two subgraphs of $G$ are said to differ by a $Z$-transformation if their
symmetric difference consists of the boundary edges of a single cell---and
if each subgraph contains an alternating set of the boundary edges of that
cell.  We say two elements of $S_{\phi}$ are adjacent if they differ by a
single $Z$-transformation.  The problem is to understand the connected
components of $S_{\phi}$.

In the case that $\phi$ is identically one (and thus $S_{\phi}$ is the set
of perfect matchings of $G$) the problems has been solved explicity for
polyomino graphs, simply connected hexagon graphs, and other families
without the aid of height functions (see \cite{L}, \cite{RF}, \cite {Z},
\cite{ZZ3} and \cite{ZGC}).

Using height functions, however, we can {\it always} determine the
structure of $S_{\phi}$ and enumerate its connected components.  The key
is to observe that two subgraphs differ by a $Z$-transformation if and
only if their corresponding height functions are equal everywhere except
on a single face where they differ by one.  In fact, we will extend our
definition to say that any two height functions differ by a
$Z$-transformation if this is the case.  Given this formulation, the
connectivity results are almost trivial.

Similarly, several papers (including those in \cite{A}, \cite{CZ}, \cite{E},
\cite{HHS}, \cite{IPS}, and \cite{ZZ1}) have used a variety of techniques to
produce Hamiltonian cycles in grid graphs without using height functions.  
However, these algorithms and characterizations do not apply in general.  For
example, \cite{IPS} solves the Hamiltonian path problem for rectangular grids;
\cite{HHS} and \cite{S} deal with rectangular grids with one or two vertices
removed; \cite{ZZ1} defines a special class of Hamiltonian grid graphs by
inductively adjoining pairs of rows; and \cite{A} gives a linear algorithm for
the Hamiltonicity of grid graphs shaped like staircases (with possibly irregular
step sizes).

Using height functions, we can go much further.  First, we specialize to
the case that $G$ is a grid graph and $\phi$ is identically two (so that
elements of $S_{\phi}$ are cycle covers of $G$).  We will see that
applying a series of $Z$-transformations to lattice square faces sometimes
allows us to join those cycles together.  This turns out to always be
possible provided there is enough ``slack" in the corresponding height
functions; our algorithm will search through the connected components of
$S_{\phi}$ until it finds one for which the corresponding height functions
do have enough ``slack."

We will always assume that $G$ is 2-connected (otherwise, it is clearly
not Hamiltonian), and thus, the faces of $G$ are simple polygons without
repeated edges.  The complexity of our algorithms will depend on $m$, the
number of faces of $G$ that are not lattice squares.  (Intuitively, these
correspond to ``holes'' in the grid graph.)  We give a (probably not
sharp) upper bound on the runtime by $O(|F|^{2+\frac{m}{2}})$, where $F$
is the set of square faces of $G$.  In particular our algorithm is
$O(|F|^2)$ for grid graphs with no non-square faces (also called {\it
polyomino graphs}).  Similar results hold for Hamiltonian paths with fixed
endpoints and for other gridlike families of graphs embedded in $n$-holed
tori.

An algorithm very similar to the one we present was published in
\cite{UL}; this paper (which applies only to polyomino graphs and gives a
somewhat weaker bound of $O(|V|^4)$) does not discuss height functions
explicitly, but it contains many definitions and constructions that are
most naturally understood in height function language.  For example, the
somewhat unintuitive definition of the distance between $2$-factors in
\cite{UL} (based on an induction using nestings of ``alternating cycles'')
is the $l_1$ distance between the corresponding height functions.

An intuitive reason we might expect the Hamiltonicity problem to be more
difficult for grids with holes is that the path could belong to any one of
the exponentially many {\it homotopy} classes of non-intersecting cycles
in the plane minus the holes.  An interesting open problem we don't solve
is whether, given a homotopy class, one can determine in polynomial time
whether there is a Hamiltonian cycle in that homotopy class.  To
compensate, we provide a related result that is neither strictly stronger
nor strictly weaker: we produce a correspondence between cycle partitions
and network flows and provide a polynomial algorithm for the existence of
a Hamiltonian cycle whose corresponding network flow lies in a particular
{\it homology class} of the plane minus the holes.

Our last topic will be the enumeration, sampling, and counting of the
Hamiltonian cycles (or Hamiltonian paths) for a given grid graph using
Markov chains.  Although the convergence rates of the Markov chains are
not known in general, we can provide heuristics and upper bounds in
special cases.

Finally, as a potential application, we point out that theoretical
chemists and statistical physicists make frequent use of Hamiltonian paths
and cycles in grid graphs to model proteins and other heteropolymers.  
These results may prove useful in analyzing the phase transitions and
other properties of two-dimensional single-chain heteropolymer models.

\section{Height Functions: Definitions and Basic Results}

The results in this section are standard material from combinatorial
topology and the cycle decomposition of planar electric flows; we include
them to fix notation and make the paper self-contained.

\subsection{Defining Chains and Homology}

From here on, we assume that $G$ is embedded in a plane or $n$-holed torus
and we will refer to the connected components of the surface of embedding
minus $G$ as {\it faces}.  One easily checks that if $G$ is connected and
embedded in a surface with genus as small as possible, all of the faces
are simply connected and have a sequence of (possibly repeated) oriented
edges $e_1, \ldots, e_k$ as a clockwise oriented boundary.  Because of the
orientability of our surface, if a face contains the same edge twice on
its boundary, it must be oriented in opposite directions each time.  If
none of the faces of $G$ has repeated edges on its boundary, then all face
boundaries are simple polygons, and we say $G$ is a {\it polygon graph}.  
We have already defined $Z$-transformations for polygon graphs, but we can
extend that definition.  If a face $f$ has repeated edges, we say two
subgraphs $H_1$ and $H_2$ differ by a $Z$-transformation at $f$ if their
symmetric difference is the set of non-repeated edges of $f$ and if one
subgraph contains all non-repeated edges oriented from black to white on
the clockwise boundary of $f$ and the other contains all such edges
oriented from white to black.

Let $S_0(G)$, $S_1(G)$, and $S_2(G)$ be the formal vector spaces generated
respectively by the vertices, oriented edges, and faces of $G$ (where in
$S_1(G)$, an edge oriented one direction is the additive inverse of the
same edge oriented the opposite direction).  For example, if $G$ has five
faces, $S_2(G)$ is the five-dimensional vector space of sums of the form
$\sum_{f \in F} a_f f$, with $a_f \in \mathbb R$.

The dual perspective is to think of an element of $S_2(G)$ as a function
from these five faces to the real numbers.  Because we will need this dual
perspective frequently, by slight abuse of notation, whenever $\tau$ is an
$2$-chain, we will write $\tau(f)$ to mean the value of $a_f$ in the sum
$\tau = \sum_{f \in F}a_f f$, (i.e., the inner product of $\tau(f)$ and
$f$).  We follow similar conventions for $S_1(G)$ and $S_0(G)$.  Elements
of $S_i(G)$ are called {\it $i$-chains}.  There is an obvious
correspondence between 1-chains and electric flows.

Next, we define a boundary operator:

$$d: S_1(G) \mapsto S_0(G)$$
$$d: S_2(G) \mapsto S_1(G)$$

as follows; if $e = (v,w)$ is an edge (oriented from vertex $v$ to vertex
$w$), then $d(e) = w - v$.  If $f$ is a face, then $d(f)$ is the sum of
the edges on the boundary of $f$, clockwise oriented.  (``Clockwise" is
well-defined because the plane and $n$-holed torus are orientable.)  Note
that if the boundary of $f$ contains a repeated edge, both that edge and
its inverse will be included in the boundary sum and will cancel each
other out.

If $\omega$ is a 1-chain, we say $\omega$ is {\it closed} if $d(\omega) =
0$.  We say $\omega$ is {\it exact} if $\omega = d(\tau)$ for some chain
$\tau$.  Since, $d^2(f)$ is clearly equal to zero for any single face $f$,
$d^2(\tau) = 0$ for any 2-chain $\tau$.  Hence, every exact one-chain is
closed.

The homology $H_1(G)$ is defined to be the vector space of closed 1-chains
mod the space of exact 1-chains.  In other words, homology represents the
failure of the sequence $$0 \mapsto S_2(G) \mapsto^d S_1(G) \mapsto^d
S_0(G) \mapsto 0$$ to be exact.  The most important result from
combinatorial topology we use is the following:

\begin{theorem}
The dimension of the homology space depends only on the topological
structure of $T$, the surface of embedding, not on the graph $G$.  The
homology is trivial if $T$ is the plane and $2n$-dimensional if $T$ is
the $n$-holed torus.
\end{theorem}

The reader who wishes to verify these facts can do so with the following
steps.

\begin{enumerate}
\item Verify that the dimension of $dS_2(G)$ is $|F|-1$ by checking that
the kernel in the space $S_2$ of the map $d$ is one-dimensional, spanned
by $\sum_{f \in F} f$.
\item Verify that the dimension of $dS_1(G)$ is $|V| - 1$ by showing that
for any two vertices, $v, w$, the 0-chain $v-w$ is in $dS_1(G)$ and then
showing that these elements span a codimension one subspace of
$S_0(G)$.
\item Conclude that the kernel of $d$ in $S_1(G)$ has dimension $|E|
- |V|+1$ and the homology has dimension $2 - |F| + |E| - |V|$.  Compute
this value using the Euler characteristic.
\end{enumerate}

\subsection{Defining $Z$-transformations and $S_{\phi}$ in Terms of Chains}

We refer to vertices of $G$ as {\it black} and {\it white} depending on
which of the two partite classes they belong to.  Given a subgraph $H$ of
$G$, we then define the 1-chain $\omega_H$ to be the
sum of all the edges in $H$, oriented from black to white.  Then $H \in
S_{\phi}$ if and only if
$$ d\omega_H = \sum_{v \in G} \phi(v)\epsilon(v)v $$
where
$$
\epsilon (v) =
\begin{cases}
1 & \text{$v$ is a black vertex} \\
-1 & \text{$v$ is a white vertex}
\end{cases}
$$

From this it follows also that if $H_1$ and $H_2$ are both in $S_{\phi}$,
$d\omega_{H_1} - d\omega_{H_2} = 0$.  In other words, $\omega_{H_1} -
\omega_{H_2}$ is always closed (though not necessarily exact, unless $G$
is planar).  Furthermore, two subgraphs $H_1$ and $H_2$ differ by a
$Z$-transformation if and only if $ \omega_{H_1} - \omega_{H_2}$ is equal
to $df$ or $-df$ for some single face $f$.

This implies that if one can move from $H_1$ to $H_2$ by a sequence of
$Z$-transformations (i.e., $H_1$ and $H_2$ are in the same connected
component of $S_{\phi}$) then $\omega_{H_1} - \omega_{H_2}$ can be written
as a sum $\sum_{f \in F} a_f df$ where $a_f$ are integers.  In particular,
$\omega_{H_1}- \omega_{H_2}$ is exact.

Whenever $\omega_{H_1} - \omega_{H_2}$ is exact, we will say that $H_1$
and $H_2$ are {\it homologically equivalent} and write $H_1 \sim H_2$;
since this is clearly an equivalence relation, it partitions $S_{\phi}$
into {\it homology classes}.

\subsection{Chains and General Network Flow Problems}

In the above formulation, a 1-chain $\omega$ is equal to
$\omega_H$ for some $H \in S_{\phi}$ if and only if
\begin{enumerate}
\item $d\omega = \sum \phi(v)\epsilon(v)v$ 
\item For each black to white oriented edge $e$, $\omega(e) \in \{0, 1\}$.
\end{enumerate}

Equivalently, we might require that $\omega$ represent an electrical
network flow such that each black vertex is a source of $\phi(v)$ units of
current, each white vertex is a sink of $\phi(v)$ units, and each edge
conducts one or zero units from black to white.  This is best understood
as a special case of the following well-known constrained network flow
problem:

Given a graph $G$, an upper bound $\delta(e)$ and a lower bound (given by
$-\delta(e^{-1})$) on the amount of current passing through each directed edge
$e$, and a function $\phi$ on the vertices dictating the net amount of current
(positive or negative) flowing into that vertex, describe the set of (real or
integer-valued) flows with this property, i.e., the set of one-chains $\omega$
such that $d\omega = \sum \phi(v) v$ and $\omega(e) \leq \delta(e)$ for each 
directed edge $e$.

The first step is to reduce to the case that $\phi$ is identically zero (and
hence our one-chains are closed) as follows.  Let $\alpha$ be any (integer or
real valued) flow such that the amount of flow into $v$ is $\phi(v)$ for each
$v$.  Instead of seeking $\omega$, we can look for $\overline{\omega} = \omega -
\alpha$ such that $d\overline{\omega} = 0$ and $\overline{\omega}(e) \leq
\overline{\delta}(e)$ for all $e$, where $\overline{\delta} = \delta - \alpha$.

Next, we treat the set of closed $\overline{\omega}$ one homology class at a
time.  (If $G$ is planar, there is only one homology class.)  Accordingly, we
restrict ourselves to $\omega$ such that $\overline{\omega}$ is exact.  
(Different choices of $\alpha$ correspond to different homology classes.)  Then
for each such $\overline{\omega}$, there is a 2-chain $\tau$ such that $d\tau =
\overline{\omega}$.  Viewed in the dual sense as a function on the faces of $G$,
$\tau$ is the {\it height function} discussed in the introduction, and it is
uniquely determined up to a constant.  To eliminate ambiguity, we will require
that $\tau(f_0) = 0$ for some reference face $f_0$.  If $f_1$ and $f_2$ share an
edge $e$ (oriented with $f_1$ on the left), the condition that
$\overline{\omega}(e) \leq \overline{\delta}(e)$ is equivalent to the condition
that $\tau(f_2) - \tau(f_1) \leq \overline{\delta}(e)$.  We have
now reduced the characterization of flows in a homology class to the
characterization of height functions satisfying difference restrictions.

To help us better understand the restrictions, for adjacent $f_1$ and
$f_2$, we define $d(f_1, f_2) = \overline{\delta}(e)$, the largest allowed
value of $\tau(f_2) - \tau(f_1)$.  Next, suppose $f_1$ and $f_2$ do not
necessarily share an edge, but $P = (f_1, f_{\sigma(1)}, \ldots,
f_{\sigma(m)}, f_2)$ is a path connecting $f_1$ and $f_2$ (so that each
$f_{\sigma(i)}$ shares an edge with the face before and after it).  Then
we define $D_P(f_1, f_2) = d(f_1, f_{\sigma(1)}) + d(f_{\sigma(1)},
f_{\sigma(2)}) + \ldots + d(f_{\sigma(m)}, f_2)$.  We then deduce the
useful bounds, $-D_P(f_2,f_1) \leq \tau_H(f_2)  - \tau(f_1)  \leq
D_P(f_1,f_2)$.  When it exists, we define $D(f_1, f_2)$ to be the minimum
value of $D_P(f_1, f_2)$ as $P$ ranges over all possible paths connecting
$f_1$ and $f_2$.  This includes single-element paths; thus, if there are
no {\it negative cycles} (i.e., no paths $P$ from a face $f$ to itself
such that $D_P(f,f) < 0$), then $D(f,f) = 0$.

\begin{theorem} There exists a function $\tau$ on the faces of $G$
satisfying the difference restrictions if and only if there exists no path
$P$ from a face $f$ to itself such that $D_P(f,f) < 0$.
\end{theorem}

The above bounds clearly imply that the latter condition is necessary.  
To prove that it is sufficient, we first note that if there are no
negative cycles, then for any $f_1$ and $f_2$, a minimum $D_P(f_1,f_2)$
must exist.  To see this, note that if $P'$ is a path connecting $f_2$
back to $f_1$ and $P+P'$ is the concatenation of $P$ and $P'$, then
$$0 \leq D_{P+P'}(f_1, f_1) = D_P(f_1, f_2) + D_{P'}(f_2, f_1)$$
If we can make $D_P(f_1, f_2)$ arbitrarily small, we obtain a
contradiction by choosing $P$ so that $D_{P+P'}(f_1, f_1)$ is negative.

Now, to prove existence of a $\tau$ satisfying our restrictions, we show
that there is one choice of $\tau$ that achieves its maximal allowed value
on every face: that is, $\tau(f) = D(f_0, f)$ for all $f$.  To see that
this function satisfies our restrictions, suppose $f_1$ and $f_2$ are
joined by an edge $e$ with $f_1$ on the left.  Then one easily checks that

$$D(f_0, f_2) \leq D(f_0, f_1) + d(f_1, f_2) = D(f_0,f_1) + 
\overline{\delta}(e)$$
$$D(f_0, f_1) \leq D(f_0, f_2) + d(f_2, f_1) = D(f_0, f_2) - \overline{l}(e)$$
and hence
$$\tau(f_2) - \tau(f_1) \leq \delta(e)$$

From this construction of $\tau$, we get the following corollary.

\begin{corollary}
If a reference face $f_0$ is required to have $\tau(f_0) = 0$, then there
exists a $\tau$ satisfying our difference restrictions such that $\tau$
achieves its maximal (or symmetrically, minimal) possible value on each
face of $G$.  In other words, if we write $\tau_1 \geq \tau_2$ whenever
$\tau_1(f) \geq \tau_2(f)$ for all $f$, the resulting poset of allowable
height functions, contains maximal and minimal elements.
\end{corollary}

The distance function $D$ can be thought of as representing a
(non-symmetric and possibly negative) cost of traveling from $f_1$ to
$f_2$ and can be efficiently computed by standard shortest path and
negative cycle detection algorithms, such as Dijkstra's algorithm (see
\cite{AMO}) and its many variants.

Clearly, every homology class of $S_1(G)$ corresponds to the class of
height functions arising from some choice of $\alpha$.  Now, fix
$\alpha_0$ to be any chain satisfying $d\alpha_0 = \sum \phi(v) v$, and
let $c_1, \ldots, c_n$ be a generating basis of the elements of the
homology space of $G$; (such a basis can be computed with simple linear
algebra, or $Z$-module reduction in the integer-valued case).  We can now
enumerate the homology classes by determining for what values of $a =
(a_1, \ldots, a_n)$ the one-chain $\alpha = \alpha_0 + \sum a_i c_i$ gives
rise to negative cycles in the corresponding restricted difference
problem.

It is clearly enough to consider non-self-intersecting paths $P$, of which
there are finitely many.  Given a path $P$ from a face $f$ to itself, note
that the value $D_P(f,f)$ has affine dependence on $a$, and thus the
restriction that $D_P(f,f) \geq 0$ forces $a$ to lie one side of a
hyperplane.  Hence, for the real (integer-valued) flow problem the
allowable $a$'s correspond to real (integer-valued) vectors inside a
convex polytope.  Furthermore, for any allowable $a$ and $a'$ and any $i$,
one easily checks that $|a_i - a'_i|$ is bounded above by $|F|M$, where
$M$ is the maximum size of a restricted difference interval.  Thus, in the
integer-valued case, $(|F|M)^n$ is an upper bound on the number of
homology classes.  In particular, for the restricted vertex degree
problem, $|F|^n$ is a (usually very weak) upper bound for the number of
homology classes of $S_{\phi}$.

\section{Connectedness Results}
\subsection{Connected Components of Height Function Spaces}

In this section, we study the connectedness of the set of flows in a given
homology class; recall that once the homology class is given, allowable
flows are in one to one correspondence with functions $\tau$ on the faces
of $G$ such that whenever $f_1$ and $f_2$ share an edge $e$ (oriented with
$f_1$ on the left), $\overline{l}(e) \leq \tau(f_2) - \tau(f_1) \leq
\overline{u}(e)$.

Now, given a set of faces and difference restrictions, we say two faces
$f_1$ and $f_2$ are in the same face cluster (or $f_1 \sim f_2$) if
$D(f_1,f_2) = -D(f_2,f_1)$.  Equivalently, $f_1 \sim f_2$ if for some path
$P$ from $f_1$ to itself which passes through $f_2$, $D_P(f_1, f_1) = 0$.  
This is clearly an equivalence relationship, and moreover, if $f_1 \sim
f_2$, one easily checks that $D(f_1,f_3) = D(f_1,f_2) + D(f_2,f_3)$ and
$D(f_3, f_1) = D(f_3,f_2) + D(f_2,f_1)$ for any third face $f_3$.  
Whenever $f_1 \sim f_2$, the difference $\tau(f_1) - \tau(f_2)$ is fixed
to be exactly $D(f_1, f_2)$ for any allowable height function; thus, it is
always impossible to apply $Z$-transformations to an allowable $\tau$ at
$f_1$ or $f_2$.

The following result is now easy to prove, and is based on the intuitive
notion from linear programming that if a convex polytope is
full-dimensional in its space of embedding, any two points (or at least
interior points) can be connected by paths that only move in one
coordinate direction at a time.

\begin{theorem} The $Z$-transformation graph on $S_{\phi}$ is connected if
and only if there is at most one face cluster containing more than one
face.  If there are $k$ face clusters ($k > 1$) with representatives
$f_0, \ldots, f_{k-1}$, then (taking $f_0$ as the reference face so
$\tau(f_0)=0$), the connected components of the $Z$-transformation graph
are in one to one correspondence, with the set of values $\tau(f_1),
\ldots, \tau(f_{k-1})$ satisfying $-D(f_j,f_i) \leq \tau(f_j) - \tau(f_i)
\leq D(f_i,f_j)$ for all $0 \leq i,j, \leq k-1$. \end{theorem}

First, since we cannot apply $Z$-transformations to make changes to
$\tau(f)$ when $f$ is on a fixed face, any $\tau_1$ and $\tau_2$ are
clearly in different connected components if they disagree on any of the
$f_1, \ldots, f_{k-1}$.

For the other direction, we will need some definitions.  For each path $P$
from a face $f$ to itself that includes at least one face other than $f$,
define the length $|P|$ to be the number of edges of the path (so $(f,
f_1, f_2, f)$ has length three).  Then let $\beta$ be the minimal value of
$D_P(f, f)/|P|$ as $f$ ranges over all faces and $P$ ranges over all
non-trivial paths from $f$ to $f$ that do not intersect themselves except
at their endpoints.  Next, define the distance $|\tau_1 - \tau_2| =
\sum_{f \in F} |\tau_1(f) - \tau_2(f)|$.

We claim that with at most $\frac{|\tau_1-\tau_2|}{\beta} + |F|$
$Z$-transformations, we can transform $\tau_1$ into $\tau_2$.  To see
this, first let $\gamma$ be the maximal value assumed by $\tau_1(f) -
\tau_2(f)$ as $f$ ranges over the faces in $F$.  We can assume without
loss of generality that this value is positive, and $f$ is a face on which
this value is achieved.  Next, let $m$ be the minimal value such that
decreasing the value of $\tau_1(f)$ produces an allowable height function.  
Then if $m \geq \tau_1(f) - \tau_2(f)$, we decrease the value of
$\tau_1(f)$ by $\tau_1(f) - \tau_2(f)$, making the two height functions
equal at that face.  If this is not the case but $m\geq\beta$, we decrease
the value of $\tau_1(f)$ by $m$, reducing the distance $|\tau_1 - \tau_2|$
by at least $\beta$.  Finally, suppose $m < \tau_1(f) - \tau_2(f)$ and $m
< \beta$.  Then there must be a neighboring $f_1$ such that decreasing
$\tau_1(f_1)$ by $m$ causes $\tau_2 - \tau_1$ to be greater than
$d(f_1,f_2)$.  In other words, there is a neighboring $f_1$ such that,
$d(f,f_1) - (\tau_1(f_1) - \tau_1(f)) < m < \beta$.  Next, we go through
the above steps again, either reducing the distance by at least $\beta$,
creating equality on one face, or finding an $f_2$ with $d(f_1,f_2) -
(\tau_1(f_2) - \tau_1(f_1)) < m \leq \beta$.  The sequence $f, f_1, f_2,
\ldots$ we produce here can never intersect itself, because if $f_k,
\ldots, f_l$ were are cycle $P$ formed by this process, we would have
$$D_P(f,f) = \sum_{i=k}^{l-1} d(f_i,f_{i+1}) - (\tau_1(f_i) -
\tau_1(f_{i+1})) = \sum_{i=k}^{l-1} d(f_i,f_{i+1}) < (l-k)\beta$$ and this
contradicts our assumption about the minimality of $\beta$.  Hence, we can
either decrease the distance by $\beta$ or produce equality on a face at
every step, and we can apply a similar procedure for the faces on which
$\tau_1(f) < \tau_2(f)$; from this we deduce the bound on the number of
steps.  Finally, to see that there exists at least one possible height
function for each allowable set of values $\tau(f_1), \ldots,
\tau(f_{k-1})$, we can use an extrapolation argument to set each $\tau(f)$
equal to its maximal possible value.  That is, $\tau(f)$ is the maximum of
$\tau(f_i) + D(f_i,f)$ as $i$ ranges from $0$ to $k-1$.  The proof that
this yields a valid height function is analogous to that given earlier for
height function existence. $\square$

In the integer-valued case, every such step decreases the distance between
$\tau_1$ and $\tau_2$ by at least one, so we get the following corollary.

\begin{corollary}
In the case of integer-valued network flows, if $\tau_1$ and $\tau_2$
agree on all of the face clusters with more than one element, $\tau_1$
and $\tau_2$ can be connected by at most $|\tau_2 - \tau_1|$
$Z$-transformations.
\end{corollary}

Furthermore, when we are interested in random sampling or Coupling from
the Past, we can always make our space of flows connected by allowing
$Z$-transformations on face clusters instead of just individual faces.  
Further modifications may be made to increase the coupling speed in many
cases, but these are beyond the scope of this paper.

\subsection{Restricted Vertex Degree Subgraph Connectedness}

In the Restricted Vertex Degree Subgraph case, we can restate the
connectedness result using the following language.  We say an edge of $G$
is {\it free} if there exist $H_1, H_2 \in S_{\phi}$ such that the $e$ is
contained in $H_1$ and not $H_2$.  Otherwise, we say the edge is {\it
fixed}; i.e., it is necessarily included (not included) in each element of
$S_{\phi}$.  Similarly, we say an edge $e$ of $G$ is {\it homology class
free} with respect to a given homology class if there exist $H_1, H_2$ in
that homology class such that $e$ is contained in $H_1$ and not $H_2$;
otherwise, $e$ is {\it homology class fixed}.  Clearly, two adjacent faces
are in the same face cluster if and only if the edge between them is
homology class fixed-and in the planar case, homology class fixed edges
and fixed edges are equivalent.

The notion of free and fixed edges is used in organic chemistry.  If we
treat carbon atoms as vertices, single or double bonds as edges in $G$,
and double bonds as edges in $H$, benzene molecules exist only when $G$ is
a hexagon graph that admits a perfect matching; when multiple matchings
are present, the double bond configuration occupies a quantum
superposition of possible matchings, and carbon atoms joined by free edges
of $G$ have different bond lengths from those joined by fixed edges.

Using this definition, we have another formulation of our connectedness
results:

\begin{theorem} The following are equivalent:
\begin{enumerate}
\item There exists a 1-chain $\alpha$ defined on $G$ such that $d\alpha =
\sum \phi(v) \epsilon(v) v$ and each edge $e$ oriented from black to
white satisfies $0 < \alpha(e) < 1$.
\item There are no fixed edges in $G$.
\end{enumerate}
\end{theorem}

To see that the first implies the second, first, as in previous sections,
we use $\alpha$ to define $\overline{\omega}_H = \omega_H - \alpha$ and we
choose $\tau_H$ to be the $2$-chain such that $\tau_H(f_0) =0$ for a
reference face $f_0$ and $d\tau_H = \overline{\omega}_H$.  Note that the
statement implies that for any pair of adjacent faces, $f_1$ and $f_2$, $0
< |d(f_1, f_2)| < 1$.  In particular, this implies that for any
non-trivial path $P$ from a face $f$ to itself, $D_P(f,f)$ is positive,
and hence each face cluster contains only one face.

For the converse, if there are no fixed edges on $G$, then let $H_1,
\ldots, H_k$ be some set of subgraphs in $S_{\phi}$ such that edge is
contained and not contained in at least one subgraph.  Then $\alpha =
\frac{1}{k}\sum \omega_{H_i}$ satisfies the requirements. $\square$

Note that if we let $H_1, \ldots, H_k$ run through all the subgraphs in
$S_{\phi}$, then $\alpha \frac{1}{k}\sum \omega_{H_i}$ and hence
$\frac{1}{k} \sum \overline{\omega}_{H_i}=0$.  We then conclude that
$\frac{1}{k}\sum \tau_{H_i} = 0$, because it has zero as its boundary and
because $\tau_{H_i}(f_0)$ is always equal to zero.  This is the unique
choice of $\alpha$ such that the average of the corresponding height
functions is zero.  Clearly, the problem of determining this $\alpha$ ---
the average value of $\omega_{H_i}$ --- is equivalent to the problem of
determining the average height function when a different, fixed $\alpha$
is chosen.

Using the same proof, we can generalize to clusters of faces.

\begin{theorem} The following are equivalent:
\begin{enumerate}
\item There exists a 1-chain $\alpha$ defined on $G$ such that $d\alpha =
\sum \phi(v) \epsilon(v) v$ and for a partition of the set of faces $F_1,
\ldots, F_k$, each edge $e$ separating members of different $F_i$
oriented from black to white satisfies $0 \leq \alpha(e) < 1$.
\item No edge separating members of different $F_i$ is fixed.
\end{enumerate}
\end{theorem}

This follows from the fact that if an edge $e$ were fixed, then we would
have $f_1 \sim f_2$ for the two opposing faces and hence there would exist
a closed path $P$ starting and ending at $f_1$ and passing through $f_2$
with $D_P(f_1, f_1) = 0$.  However, we know that all the terms in the sum
for $D_P(f_1,f_1)$ are equal to or greater than zero and at least one (one
separating an element of the first $F_i$ from the second) is positive.
$\square$

Similarly, we can state the result for particular homology classes.  This
latter formulation is useful because the second condition implies that all
subgraphs of $G$ corresponding to the homology class defined by $\alpha$
(i.e., all $H$ such that $\omega_H - \alpha_0$ is exact) can be connected
by a series of $Z$-transformations.

\begin{theorem} The following are equivalent:
\begin{enumerate}
\item There exists a 1-chain $\alpha$ defined on $G$ such that $d\alpha =
\sum \phi(v) \epsilon(v) v$, $\alpha-\alpha_0$ is exact, and for a
partition of the set of faces $F_1,\ldots,F_k$, each edge $e$ separating
members of different $F_i$ oriented from black to white satisfies $0 \leq
\alpha(e) < 1$.
\item No edge outside of $\mathcal E$ is homology class fixed with
respect to the homology class of $H$ for which $\omega_H  - \alpha_0$ is
exact.  In other words, no faces from different $F_i$ are part of the
same face cluster.
\end{enumerate}
\end{theorem}

Finally, if $f_0$ is the unbounded face of a planar graph $G$, we will
sometimes be given a natural choice of $\alpha$ such that $\alpha -
\omega_{H_0}$ is closed at all interior vertices of $G$ but not
necessarily at those vertices $v_1, \ldots, v_k$ which border on $f_0$; in
this case, our height functions cannot be well defined on $f_0$.  
However, suppose we add a new vertex $v$ to the middle of $f_0$ together
with an edge from $v$ to each vertex on the boundary of $f_0$, thus
dividing $f_0$ into triangles, one triangle for each edge.  Then for each
edge $e_i$ connecting $v$ to $v_i$, we define $\alpha(e_i)$ so that
$\alpha-\omega_{H_0}$ is closed at $v_i$.  Having done this, it follows
that $\alpha$ is also closed at $v$ (since $v$ cannot be a source or sink
if no other vertex is a source or a sink), and our height functions
$\tau_H$ are now well-defined for any $H$ in $S_{\phi}$; and if we choose
one of the triangles in the division of $f_0$ as a reference face, we see
that the value of $\tau_H$ on each of these triangles is independent of
the choice of $H$.  We will refer to these as the {\it boundary values},
noting that one such value corresponds to each edge of $f_0$.

One reason for adopting the above perspective is to force $\alpha$ to
satisfy some regularity conditions.  For example, suppose $G$ is a
bipartite planar graph such that every vertex not on that unbounded face
has exactly has $k$ edges and that $\phi$ is identically equal to $j$
(with $0 < j < k$).  Then can we set $\alpha$ identically equal to $j/k$
on each edge oriented from black to white.  For example, in the case of
perfect matchings on a polyomino graph, we can take $\alpha(e)$ to be
$1/4$ everywhere for each edge oriented from black to white.  This is
essentially the approach used to describe height functions in, for
example, \cite{CEP}.  Similarly, if $G$ is a hexagonal graph with no
``holes" (i.e., all unbounded faces are hexagons), we can understand the
perfect matchings of $G$ by taking $\alpha = 1/3$.  Applying the
connectedness arguments of the previous section, we can also deduce the
following:

\begin{theorem}
Let $G$ be a bipartite planar graph such that every vertex not on the
unbounded face has exactly $k$ edges and $\phi$ is constant.  Then
$S_{\phi}$ is connected.
\end{theorem}

In the case of partitions of polyomino graphs into cycles, we can take
$\alpha(e)$ to be identically $1/2$; thus, the partitions of the graph
into cycles correspond to the ways of choosing $\tau$ on the interior
faces such that whenever $f_1$ and $f_2$ are neighbors,
$|\tau(f_1)-\tau(f_2)|=1/2$; and if $f$ shares an edge with the unbounded
face, $f$ differs by plus or minus $1/2$ from the boundary value of that
edge.  In our grid graph illustrations in the next section, we multiply
$\tau$ by two so that we can deal with integer boundary values and
increments of one.

\section{Computing the Hamiltonicity of Grid Graphs}

\subsection{Fixed Faces and Grid Graphs}

In this section, we let $G$ be a {\it grid graph}, i.e., an induced
subgraph of $\mathbb Z^2$, and take $\phi$ to be identically $2$.  Thus,
$S_{\phi}$ consists of the spanning subgraphs $H$ of $G$ that partition
$G$ into some number of cycles.  We assume $G$ is a polygon graph and
define a {\it hole} in the graph $G$ to be a bounded face of $G$ with more
than four edges.  (All two-connected grid graphs are polygon graphs, and
graphs that fail to be two-connected cannot be Hamiltonian, so we lose no
generality with our assumption.)

The first question we ask is, what are the connected components of
$S_{\phi}$ if we only allow $Z$-transformations on the faces of $G$ that
are squares?  Suppose $H_1$ and $H_2$ are two partitions of $G$ into
cycles.  Then $\omega_{H_1} - \omega_{H_2}$ is closed and hence exact.  
Choose $\tau$ so that $d\tau=\omega_{H_1}-\omega_{H_2}$ and $\tau(f_0) =
0$, where $f_0$ is the unbounded face.  A clear necessary condition for
$H_1$ and $H_2$ to be connectable by $Z$-transformations on squares is
that $\tau$ be equal to zero on all of the holes of $G$.  We claim that
this is also sufficient.

\begin{theorem}
Using the above notation, $H_1$ and $H_2$ are in the same component of
$S_{\phi}$ if and only if $\tau(f)=0$ for all holes $f$.
\end{theorem}

Suppose without loss of generality that $\tau(f) > 0$ for some
square $f$, and let $F_+$ be the set of all faces of $G$ on which $\tau$
achieves its maximal value; note that by assumption, all of these faces
are squares, not holes.

Recall that if two faces $f_1$ and $f_2$ share an edge $e$, ($f_1$ on the
left with $e$ is oriented black to white), then $\omega_{H}(e)$ can take
on two values depending on whether $e \in H$.  We say $f_2$ is {\it higher
(lower)} than $f_1$ (with respect to $H$), or $f_1 \prec f_2$, if
$\omega_{H}(e)$ takes on the larger (smaller) of these two values.  We
will also say that $f$ is a {\it local maximum (minimum)} if it is higher
(lower) than all of its neighbors.  Clearly, one can apply a $Z$
transformation to $H$ at $f$ if and only if $f$ is a local maximum or
minimum with respect to $H$.

We will now show that we can always find a face $f_0 \in F_+$ such that
$f$ is higher than all of its neighbors, and hence, by applying a
$Z$-transformation to $H_1$ at $f_0$, we can decrease $|\tau|$ by one.  
First, define $h(f)$ to be the number of faces adjacent to $f$ that are
higher than $f$ minus the number that are lower.  Because the
contributions from edges between squares of $F_+$ cancel out, $\sum_{f \in
F_+} h(f)$ is the number of edges of squares in $F_+$ that border on faces
not in $F_+$.  In particular, the sum is positive, and the average value
of $h(f)$ is greater than zero.  And this implies that there exists at
least some $f_0$ in $F_+$ which is higher than at least three of its
neighbors---say, those on the south, west, and east.  If $f_0$ is also
higher than its neighbor on the north, we are done.  If not, then let
$f_1$, be the face immediately north of $f_0$; because it is higher than
$f_0$, it is necessarily in $F_+$ and hence necessarily a square.

Like $f_0$, $f_1$ is higher than its neighbor to the south.  We claim that
$f_1$ is also higher than its neighbors on the east and west.  If, say,
the western neighbor were higher, there would be an increasing sequence
$f_w,f_0,f,f'$ of faces all incident to a single vertex $v$; this would
require either that all three of edges dividing these faces be in $H$ or
all three not in $H$.  This is impossible since $\phi(v)=2$.

Continuing, we choose $f_2, f_3, \ldots$ until we find a face $f_k$ in
$F_+$ that is higher than all of its neighbors.  Because the sequence
$f_1, f_2, \ldots$ must leave $F_+$ eventually, we know we find such an
$f_k$.  $\square$

Using the terminology of the previous section, this theorem is equivalent
to the statement that are no face clusters of $G$ that consist entirely of
squares.

Thus, if $G$ has $m$ holes, then when we only allow $Z$-transformations on
squares, the connected components of $S_{\phi}$ are indexed by the integer
lattice points in a convex polyhedron of dimension at most $m$.  How many
of these components are there?  Well, if we take $F$ to be the set of
square faces of $G$, then there is a path from any hole to the boundary
face or some other hole that passes through at most $\sqrt{|F|}$ faces.  
We say two holes (or a hole and the boundary face) are ``close" if they
are adjacent or if there is a path from one to the other that passes
through at most $\sqrt{|F|}$ squares.  One can easily show that the set of
holes and the unbounded face are connected under the closeness relation.  
The number of values that $\tau_{H}$ can take on any given hole once the
value at a close hole is determined is bounded by $2 \sqrt{|F|}$, and it
follows that there are at most $2^m|F|^{m/2}$ connected components of
$S_{\phi}$.

In the remaining subsections, we will only deal with $Z$-transformations
on squares, not holes, and we use only these transformations to describe
adjacencies between members of $S_{\phi}$.  Our approach to determining
the Hamiltonicity of $G$ will be to search through each connected
component of $S_{\phi}$ for a Hamiltonian path.  Let $p(H)$ be the number
of cycles into which $H$ partitions $G$.  In this section, we will write
$H \sim H_0$ when $H$ and $H_0$ are in the same component of $S_{\phi}$.  
(Equivalently, $H \sim H_0$ when $\overline{\omega}_H$ and
$\overline{\omega}_{H_0}$ are homologically equivalent in the plane minus
the holes.)

In coming sections, we assume that we begin with an $H_0 \in S_{\phi}$ and
will show that in $O(|F|^2)$ steps, we can find the $H\sim H_0$ for which
$p(H)$ is minimal.  It will follow that in $O(|F|^{m/2 + 2})$ steps, we
can find an $H$ for which $p(H)$ is minimal overall---in particular, we
will know whether $G$ is Hamiltonian.

\subsection{Height Functions for Illustrations}

If $G$ doesn't have any holes, then instead of considering the closed
1-chain $\omega_H - \alpha$ and the $\tau$ that has this as its boundary,
it will be more instructive for us to consider the 1-chain
$\overline{\omega}_H = 2\omega_H-\beta$, where $\beta$ is the sum of all
edges in $G$ oriented from black to white.  Thus, if $e$ is oriented black
to white, $\omega(e) = 1$ if $e$ is in $H$ and $-1$ otherwise.

Since $d\overline{\omega}(v)$ is not necessarily zero when $v$ is on the
unbounded face, we must adopt the ``boundary value" approach described
earlier; that is, we compute the boundary value for all the edges of the
unbounded face.  (If two boundary edges share a vertex $v$, the difference
between the value of any $\tau_H$ on the two edges is determined from the
fact that $\phi(v)=2$.)  Now the set of allowable $\tau_H$ is the set of
functions on the other faces of $G$ whose value on each face differs from
the value of its neighboring faces (or boundary edges) by $1$ or $-1$.  
The mnemonic to keep in mind is that as we move from square to square with
a black vertex on our left, $\tau_H$ goes up if $e$ is not in $H$ and down
if $e$ is in $H$.

When $G$ has holes, the problem of producing $\tau_H$ with
$\overline{\omega}_H$ as its boundary is analogous to the problem of
integrating a complex function in a multiply connected region.  Even in
this case, however, we can define $\tau_H$ on polyomino subgraphs of $G$
and treat it as a locally definable multivalued potential function.  We
will label squares of most of the examples to come with the integer values
of a locally defined $\tau_H$.  This illustrative approach will make it
easy to see when squares are locally maximal or minimal and when they are
higher or lower than their neighbors.  Note also that applying a
$Z$-transformation at a square $f$ corresponds to changing $\tau_H(f)$ by
$2$ or $-2$.

\subsection{Joining Cycles with $Z$-transformatons}

Let $B(H)$ be the set of {\it boundary faces}, that is, the faces of $G$
which contain vertices from two or more distinct cycles.  (These are not
to be confused with the primarily illustrative notion of {\it boundary
edges}.)  Now, if $f$ is a locally maximal or minimal square, let $H'$ be
the graph obtained from $H$ by applying a $Z$-transformation at $f$.  
Then if $f \in B(H)$, applying the $Z$-transformation at $f$ joins the
cycles on opposite sides of $f$ together; in this case $p(H') = p(H)-1$
and $f \not \in B(H')$.  Similarly, if $f \not \in B(H')$, applying at a
$Z$-transformation at $f$ divides the cycle with vertices on $f$ into two
cycles; in this case, $f \in B(H')$ and $p(H') = p(H)+1$.  The reader may
check this fact on the local maxima and minima in the following example.

\begin{figure} [htbp]
\begin{center}
	\includegraphics[width=8cm] {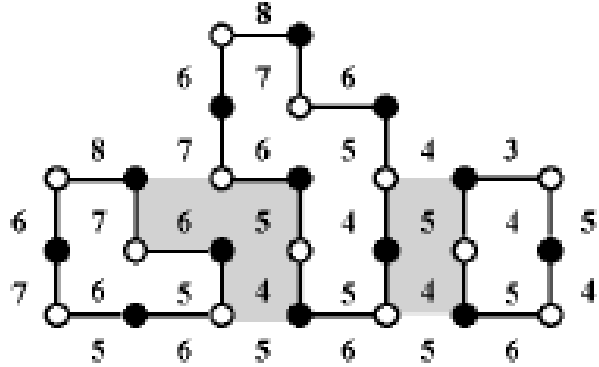}
 	\caption{Applying a $Z$-transformation at a boundary
(non-boundary) square decreases (increases) the number of disjoint cycles by one.}
\end{center}
\end{figure}
	
Since we seek the $H$ for which $p(H)$ is minimal, the obvious first step
is to hill climb until we either find an $H$ for which $p(H)=1$ (and we
are finished), or we find an $H$ such that $B(H)$ contains no locally
maximal or minimal squares.  What can we say about $H$ in the latter case?

Suppose $B(H)$ contains no locally maximal or minimal squares, and
consider a face $f \in B(H)$.  We say a vertex $v$ on $f$ is {\it
isolated} if neither of the vertices of $f$ adjacent to $v$ is in the same
cycle as $v$.  If two isolated vertices $v_1$ and $v_2$ that are incident
to $f$ are adjacent, then they belong to different cycles, and the face
which shares the edge $(v_1,v_2)$ with $f$ is necessarily a locally
maximal or minimal square in $B(H)$, a contradiction.  If $f$ contains
vertices from three or four cycles, it must contain an adjacent pair of
isolated vertices; hence, $f$ is on the boundary between exactly two
cycles.

Suppose $f$ contains one pair of vertices from each cycle.  If $v_1$ and
$v_2$ are from the same cycle, then $v_1$ and $v_2$ must be adjacent
vertices.  (Otherwise, they would be isolated.)  Similarly, the other two
vertices $v_3$ and $v_4$ are adjacent.  If both of the edges $(v_1, v_2)$
and $(v_3,v_4)$ were contained in $H$, $f$ would be locally maximal or
minimal.  If neither were contained in $f$, the two faces adjacent to $f$
containing vertices from both cycles would be locally maximal or minimal.  
Hence, $H$ contains exactly one of $(v_1,v_2)$ and $(v_3,v_4)$.  Squares
of this type are called {\it critical boundary squares}.  Inspection shows
that they are adjacent to exactly two boundary squares (or holes), and
that they are either higher than both or lower than both boundary
neighbors.  Hence, we say they are maximal or minimal {\it along the
boundary}.

Similarly, suppose $f$ has one vertex $v_1$ from one cycle and three
($v_2$, $v_3$, and $v_4$) from the other, with $v_2$ and $v_4$ adjacent to
$v_1$.  Since boundary squares adjacent to $f$ cannot contain edges from
both cycles without being maximal or minimal, one easily checks that
$(v_2,v_3)$ and $(v_3,v_4)$ must be in $H$.  Squares of this type are
called {\it corner boundary squares}, and they are not maximal or minimal
along the boundary.  We summarize this information in the following lemma.

\begin{lemma}
If $B(H)$ contains no maximal or minimal squares, then every square in
$B(H)$ is either a critical boundary square or a corner boundary square,
as illustrated in the following diagram.
\label{lemma:isolated}
\end{lemma}

\begin{figure} [htbp]
\begin{center}
	\includegraphics[width=6cm] {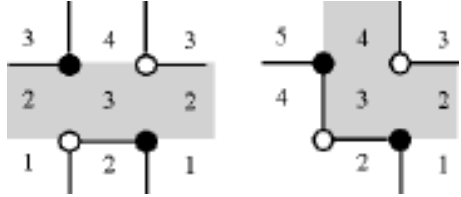}
	\caption{Critical and corner boundary squares.}
\end{center}
\end{figure}

The examples show two graphs with multiple cycles for which $|S_{\phi}| =
1$, and hence joining the cycles is impossible.  This is because each
square in these graphs is on an increasing path from one boundary edge to
another; intuitively, the height functions are pulled taut by the boundary
values, and there is not enough slack to join the separate cycles
together.

\begin{figure} [htbp]
\begin{center}
	\includegraphics[width=8cm] {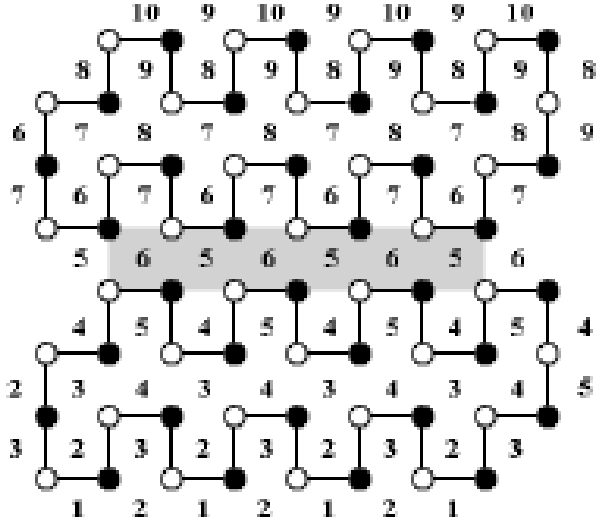}
	\caption{$|S_{\phi}| = 1$, $|B(H)|=6$, and $p(H)=2$}
\end{center}
\end{figure}

\begin{figure} [htbp]
\begin{center}
	\includegraphics[width=8cm] {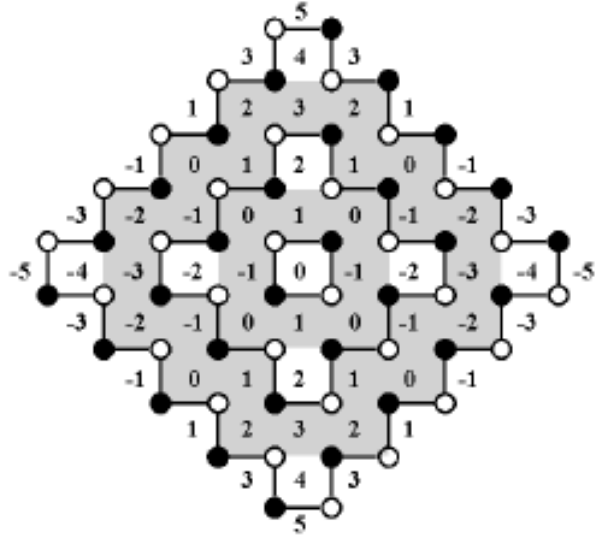}
	\caption{$|S_{\phi}| = 1$, $|B(H)|=32$, and $p(H)=3$}
\end{center}
\end{figure}

\subsection{Height Values on Boundaries}

We begin this section with a simple observation about cycle boundaries:

\begin{lemma} If the function $\tau_{H'}$ has the same value as $\tau_H$
on every boundary face in $B(H)$, two vertices that are on separate cycles
in $H$ are necessarily on separate cycles in $H'$.
\end{lemma}
Otherwise, there would have to be an edge in $H'$ connecting a pair of
vertices which belonged to different cycles in $H$.  This cannot be,
because the faces $f_1$ and $f_2$ on opposite sides of that edge are in
$B(H)$, but the difference in values on the two faces is not the same for
$\tau_H$ and $\tau_{H'}$.  $\square$

Another simple observation is the following:
\begin{lemma} If $\tau_{H'}$ has the same value as $\tau_H$ on every hole
and every critical boundary square, it necessarily has the same value on
every face in $B_H$.
\end{lemma}

To see this, note that every corner square is adjacent to two boundary
squares or holes with which it doesn't share an edge in $H$; one of these
is higher and one is lower; if these are also corner squares, they are in
turn adjacent to respectively higher or lower boundary squares or holes.  
The sequence of increasingly higher corner squares, will step
alternatively in two directions (say north and east) while the downward
sequence alternates between the other two (say south and west).  Clearly
neither sequence ever returns to a previous square in the path; hence the
upwards (downward) sequence must eventually reach a locally maximal
(minimal) boundary square or a hole.  Thus, each corner boundary square is
on an increasing or decreasing path between holes or squares that are
minimal or maximal along the boundary. $\square$

From these two lemmas, we deduce that if we hope to apply a sequence of
$Z$-transformations to $H$ to produce an $H'$ such that $p(H') < p(H)$, we
will at some point have to apply a $Z$-transformation that lowers (raises)
at least one of the critical boundary values of $H$ which is currently
maximal (minimal) along the boundary.  The next section describes a way to
do just that.

\subsection{Maximal and Minimal Rows}

Suppose that $f_0$ is a critical boundary square and that its edge $e$
containing vertices in one cycle is in $H$ and its edge $e'$ containing
vertices in the other cycle is not.  Suppose that $e$ has $f$ on its left
when oriented from black to white.  (The other case---as should be
understood thought this section---is symmetric.)  Then $f_0$ is higher
than all of its neighbor faces except for $f_1$, the other face incident
to $e'$.  Let $f_1, f_2, \ldots$ be the sequence of faces extending away
from $f$ in the direction of the edge $e'$.  Let $r$ be the smallest
integer such either $f_r$ is no longer a square or $f_r$ is lower than
$f_{r-1}$.  Clearly, $r \geq 2$.

If $f_r$ is a boundary square or not a square at all and $f_r$ is higher
than $f_{r-1}$, then we say that $f$ is a {\it stuck critical boundary
square}.  In this case, there is no way we can lower $f_0$ without first
lowering $f_r$, which is either another boundary square or a hole.

Otherwise, the sequence $f_1, \ldots, f_{r-1}$ consists of squares all of
whose vertices belong to the same cycle.  Also, $f_0, f_1, \ldots,
f_{r-1}$ is a {\it maximal square row of length r}, that is, a sequence of
$r$ consecutive squares in a row or column such that $f_i \prec f_{i+1}$
for $0 \leq i \leq r-2$ and each square in the sequence is higher than all
of its neighbors not in the sequence.  We can {\it lower} the maximal
square row by decreasing the value of $\tau_H$ by two on all elements in
the row --- or equivalently by applying our $Z$-transformations
successively to $f_{r-1}, f_{r-2}, \ldots, f_1, f_0$ --- to produce a new
subgraph $H'$.

The effect of this lowering is different depending on whether $r$ is even
or odd, as shown in the following figures.  If $r$ is odd, this has the
effect of joining the cycles separated by $f_0$; in this case $p(H') =
p(H)-1$, and $B(H')$ is a proper subset of $B(H)$.  Moreover, $\tau_H$ and
$\tau_{H'}$ agree an all elements of $B(H')$.

\begin{figure} [htbp]
\begin{center}
	\includegraphics[width=6cm] {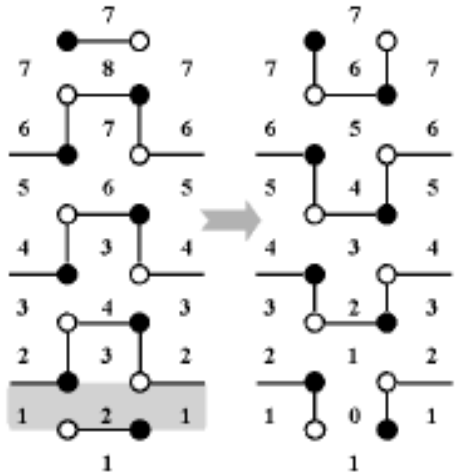}
	\caption{When a boundary square is the first square in an odd
length maximal vertex row, lowering that row joins the cycles on opposite
sides of the square.}
\end{center}
\end{figure}

If $r$ is even, the process also joins the cycles on opposite sides of
$f$.  However, by assumption, the maximal square at the top of the maximal
square row is not on a boundary.  (If it were, we could easily join two
cycles by lowering simply that square.)  Hence, lowering that square
creates an additional cycle.  In this case, $p(H') = P(H)$, and neither
$B(H')$ nor $B(H)$ is a proper subset of the other.  However, it is still
the case that $\tau(H')$ and $\tau(H)$ agree on the intersection of $B(H)$
and $B(H')$.

\begin{figure} [htbp]
\begin{center}
	\includegraphics[width=6cm] {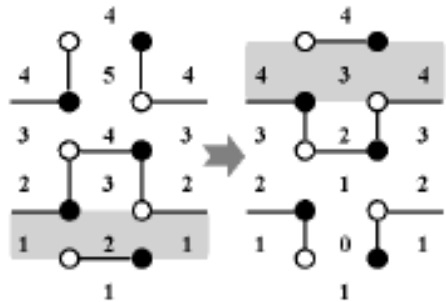}
	\caption{When a boundary square is the first square in an even
length maximal vertex row, and the last square in the row is not on a
boundary, lowering the row joins the cycles on opposite sides of the
first square and divides the cycle with vertices on the last square into
two pieces.}
\end{center}
\end{figure}

Does lowering an even length maximal row get us anywhere?  Maybe.  It
doesn't change $p(H)$, but it may be the case that it creates a new
boundary square that is the beginning of an odd length maximal or minimal
row (possibly of length one).

\subsection{Bridges}

In the scenario of the last section, suppose $r$ is even and $f_0, \ldots,
f_{r-1}$ is a maximal vertex row beginning at a boundary square $f_0$.  
(By construction, $f_{r-1}$ is not on a boundary.)  Then when we lower the
vertex $v_{r-1}$, this partitions the cycle with vertices on $f_{r-1}$
into two cycles.  At this point, $f_{r-1}$ is a minimal square on a
boundary between two cycles and $f_0, \ldots, f_{r-2}$ is an odd length
maximal square row starting on a boundary between another pair of cycles.  
I refer to a row of squares $f_0, \ldots, f_{r-1}$ of this type as a {\it
bridge}; that is, I say there is a bridge of length $r$ from $f_0$ to
$f_{r-1}$ if $f_0$ and $f_{r-1}$ are boundary squares and the squares
$f_0, f_1, \ldots, f_{r-1}$ are lined up in a row so that $f_0, \ldots,
f_{r-2}$ is an odd length maximal row of squares and $f_{r-1}$ is a
minimal vertex row (of length one).  If I raise $f_{r-1}$, thus joining
the cycles on opposite sides of that boundary point, then $f_0, f_1,
\ldots, f_{r-2}$ is no longer a maximal vertex row.  If I lower the
maximal vertex row $f_0, f_1, \ldots, f_{r-2}$, thus joining cycles on
opposite sides of the boundary point $f_0$, then $f_{r-1}$ is no longer
minimal.  Thus, I can think of a bridge as an option to join cycles at one
square or another but not both.

\begin{figure} [htbp]
\begin{center}
	\includegraphics[width=8cm] {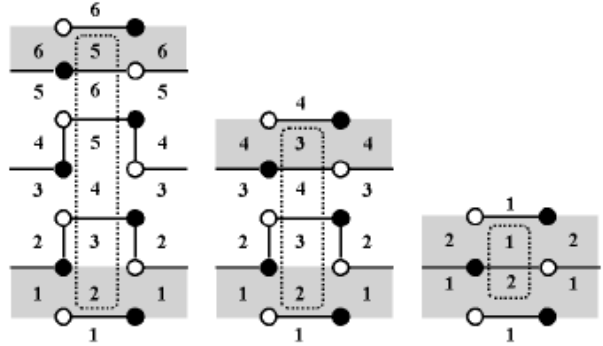}
	\caption{Bridges of length six, four, and two.}
\end{center}
\end{figure}

After constructing a bridge, I refer to the squares $f_0$ and $f_{r-1}$ as
{\it bridge end squares} and the squares $f_1,\ldots, f_{r-2}$ as {\it
inner bridge squares}.  Creating and using these bridges will be central
to our algorithm.

\subsection{The Algorithm}

This section describes a method for reducing the number of cycles of $H$
by one; we will show that this method fails only when $H$ has the minimum
number of cycles among members of its connected component of $S_{\phi}$.

We assume $H$ contains at least two disjoint cycles, and by the {\it
boundary between two cycles}, we will mean the set of all faces that are
incident to both of those two cycles.  We begin the algorithm with $B_1,
\ldots, B_k$, all of the nontrivial boundaries between pairs of cycles.  
While we are ultimately searching for an odd length maximal or minimal row
that will allow us to get rid of one of these boundaries, during the
course of this algorithm, we may temporarily create new boundaries
$B_{k+1}, \ldots, B_m$.  When a given $B_i$ is created, it represents the
set of squares between exactly one pair of cycles; as the algorithm
progresses, these may later be separated into subcycles, but that will not
concern us. We need only remember that each $B_i$ was a boundary between
exactly two cycles {\it when it was created}.  Also, each time I create a
new boundary, I will also create a bridge, as described in the previous
section, from a square on a previously existing cycle boundary to a square
on the new boundary.

The algorithm works as follows.  We start successively searching all of
the critical boundary squares $f_0$ of $H$.  If $f_0$ is stuck, or $f$ is
a square on a previously created bridge, we do nothing.  Otherwise, $f_0$
must be the beginning of a maximal (minimal) row of squares of length $r$.  
If $r$ is odd, then we lower $f_{r-1}$, thus producing a new boundary
containing $f_{r-1}$ and a new bridge from $f_0$ to $f_{r-1}$.  We say
that the new cycle boundary (the one containing $f_{r-1}$) is {\it
subordinate} to the cycle boundary containing $f_0$.  We continue this
process until one of two things happens.

\begin{enumerate}
 \item We find a critical boundary square $f_0$ (perhaps on one of the new
cycle  boundaries) that is the beginning of a maximal or minimal row of
even length.
\item Every critical cycle boundary square is either a stuck critical
boundary square or a bridge square.
\end{enumerate}

Ultimately, we'd like to show that if the first case occurs, we can lower
or raise vertex squares in a way that decreases $p(H)$ to one below its
original value.  To make this work (and simplify our proofs), we will need
two additional rules about the order in which we select our critical
squares.  Since we'd prefer not to deal with squares that are boundaries
between more than two cycles, if we ever produce a new boundary that
contains an $f$ on the boundary between three cycles, then the next square
we check in our algorithm is $f'$, the square adjacent to $f$ across an
edge connecting two isolated vertices of $f$.  We saw earlier that this
square is necessarily maximal or minimal; one can easily check that it
does not belong to any existing bridges.  This fact means that we will not
have to deal with boundary squares between three or more cycles at any
other point in our algorithm.  Except for this one case, each time a new
$B_i$ is created, if one of its squares is adjacent to a square in another
$B_i$, the two are separated by an edge in $H$.  There is one other rule
that we will add when it is motivated later on.

During this process, we will have organized our new boundaries into trees
-- under the subordination relation -- rooted at the original $B_1,
\ldots, B_k$.  Each time we add a new boundary, we add it as a leaf to an
existing node.  We will now show that once we find a cycle boundary square
on an even length maximal (minimal) row, we can join cycles in such a way
as to eliminate all of the new boundaries $B_{j+1}, \ldots, B_m$ we have
created {\it and} the boundary $B_l$ which is the root of the boundary
tree containing $B_k$, thus producing a partition of $G$ into cycles with
one fewer cycles than our original partition.

We use the following diagram to give an overview of this process. The
roots of the trees $B_1, \ldots, B_k$ represent the original boundaries
between pairs of cycles.  The $B_{k+1}, \ldots, B_m$ represent new
boundaries that are erected (in that order) during the algorithm.  Recall
that a bridge is made up of a minimal and a maximal odd length row; and we
can lower or raise exactly one of these rows to join cycles across the
boundary at one of the bridge's endpoints.  Each edge in the diagram
represents a bridge across one of the two boundaries it is connected to
--- except for the loose edge off of $B_m$ which represents a maximal
(minimal) row of squares of odd length starting at $f$ that we will call a
{\it half bridge}.

\begin{figure} [htb]
\begin{center}
	\includegraphics[width=7cm] {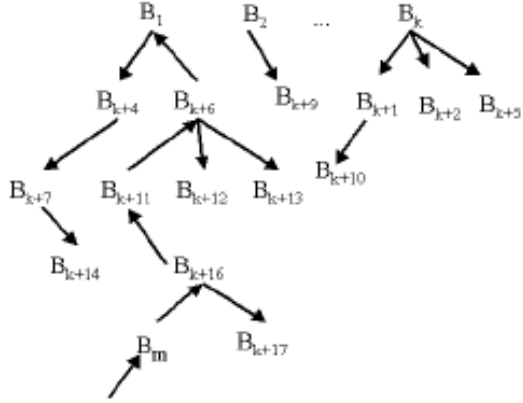}
	\caption{Using bridges and one half bridge to eliminate boundaries.}
\end{center}
\end{figure}

The arrows on the bridges point to the boundary across which we will use
the bridge to join cycles.  As shown, at each bridge we will join cycles
across the subordinate boundary except for the bridges on the path from a
root boundary to the half bridge; as shown in the diagram, we will then
have joined cycles across all of the newly created boundaries and the root
boundary $B_j$.  In order to make this precise, we will need the following
lemma, which we will prove later.  Let $M_i$ be the even length maximal
(minimal) row of squares which is lowered (raised) in order to join cycles
across the boundary $B_i$.

\begin{lemma} As sets of squares, all of the $M_i$ used by our algorithm
are disjoint.  Furthermore, no square in any one of the $M_i$ is adjacent
to any square in another $M_i$.
\end{lemma}

The lemma shows that lowering or raising one of the $M_i$ will not affect
the minimality or maximality of any of the other $M_i$.  For the remainder
of the proof, it is useful to imagine that we join cycles across the
boundaries one at a time in the reverse order of the order in which we
created them, that is, starting with $B_m$, progressing to $B_{k+1}$.

When $B_m$ was created, it separated a particular cycle into two pieces;
by joining the two cycles on opposite sides of $B_m$, we reduce $B(H)$ to
exactly what it was before $B_m$ was created.  We do the same thing by
joining cycles on opposite sides of $B_{m-1}, B_{m-2}, \ldots B_{k+1}$.  
To see this, observe that the bridge that we use to join cycles across a
given $B_i$ is either the bridge used to create that $B_i$ or a bridge
produced later in the algorithm, it follows that the bridge was created
while all $B_n$, $n < i$ were in place --- so by our construction, using
it to join cycles across $B_i$ cannot alter the value of $\tau_H$ on any
other of the $B_n$ and thus will not alter $B(H)$ other than by
eliminating the boundary squares of $B_i$.
 
Hence, after joining cycles across $B_{k+1}$, $B(H)$ is the same as it
originally was.  We then join cycles across $B_l$ to reduce the number of
cycles in $\mathcal E$ by one. $\square$

{\bf Proof of Lemma} The proof of this lemma is slightly tedious, but not
difficult.  It involves checking a number of cases.  First, we will show
that if a bridge $f_0, \ldots, f_{r-1}$ is oriented in the vertical
direction, none of the other bridges or half bridge can contain any of the
squares adjacent to $f_1, \ldots, f_{r-1}$ on the left or right.  (Of
course, the horizontal case is equivalent.)

We prove this by showing that each time we construct a new bridge or half
bridge, it cannot contain squares from any other bridge, and moreover,
that if an element of one bridge is adjacent to that of another bridge or
half bridge, the elements belong to maximal (or minimal) rows that will
not both be used (i.e., lowered or raised) by our algorithm.

So now we look at what happens when we construct a new bridge or half
bridge.  The algorithm requires us to find a cycle boundary square $f_0$
that is part of a maximal row $f_0, f_1, \ldots, f_{r-1}$, and we then
either lower $f_{r-1}$ (in the even case) to produce a new boundary and a
bridge or the entire row (in the odd case).  Now, recall that before we do
this, none of the $f_1, \ldots, f_{r-1}$ is on a boundary, and it is thus
possible to lower $f_{r-1}$ through $f_1$ without changing the value of
$\tau_H$ on any boundary square.

Now, suppose $g, g_1, \ldots, g_{k-1}$ is an existing bridge oriented from
top to bottom as in the following diagram.  As can be seen from the
diagrams, the edges separating any squares immediately to the right or
left of the $g, g_1, \ldots, g_{k-1}$ are alternatively contained and not
contained in $H$, according to the parity of their endpoints along the
bridge.  It follows that these squares are all on directed paths from one
boundary square or hole to another, and we thus cannot apply a
$Z$-transformation to any of these vertices without first changing the
value of $\tau_H$ on a boundary square.  Thus, none of these will be among
the $f_1, \ldots, f_{r-1}$.

Hence, the only way our new bridge or half bridge can contain a square
adjacent to one of the $g, \ldots, g_{k-1}$ on the left or the right is if
that square is the starting square, $f_0$.  Then $f_0$ would have to be
adjacent to either $g$ or $g_{k-1}$, since otherwise it would be in the
middle of a directed path from one cycle boundary square to another, and
hence, if critical, it would have to be stuck.  The new bridge then has to
be oriented vertically in this case since it begins on a boundary square
with another boundary square as a horizontal neighbor.  The two bridges
are adjacent only at that one square.  It is easy to check from the
diagram that in this case $f$ and $g$ must be on the same $B_k$, because
they necessarily separate the same pair of cycles.  (This follows from the
fact that there are no triple boundary squares; that case is dealt with
separately.)  However, only one of these squares (say $f$) can be part of
$M_k$; the other bridge must be used to join cycles across a different
boundary, say $B_i$, and $g$ will not be part of $M_i$.

Now, if the new bridge or half bridge is oriented horizontally, we have
shown that none of $f_0, \ldots f_{r-1}$ can be adjacent to a vertical
bridge on the left or right, and hence none of $f_0, \ldots, f_{r-1}$ can
overlap with a square in the vertical bridge.  Similarly, if the new
bridge is oriented vertically, $f_0$ cannot (by construction) be one of
$g, \ldots, g_{k-1}$, and neither of $g$ and $g_{k-1}$ can be among the
$f_0, \ldots, f_{m-1}$ (by construction); hence, the squares in the
bridges cannot overlap in this case either.

We still have one case to deal with, namely, that the new bridge $f_1,
\ldots, f_{m-1}$ may contain either the vertex $g_k$ immediately above
$g_{k-1}$ or the vertex $g_{-1}$ immediately below $g$.  First we rule out
the horizontal case.  One can easily check that both before and after a
horizontal bridge is created, all of the squares that belong to that
bridge contain at least one pair of vertices from which vertical edges
extend away from the square in either direction.  If either $g_{-1}$ or
$g_{k}$ were to be part of a bridge oriented left to right (or any part of
a half bridge except the last square), there would have to first be edges
extending upward and downward out of those squares, which cannot happen
since those edges would have to cross boundaries;  these squares cannot be
the last element of a half bridge oriented left to right because they
border an incident edge in $H$ either above or below.

Since the new bridge or half bridge cannot overlap with the $g, \ldots,
g_{k-1}$, the only remaining possibilities are that one of $g_k$ or
$g_{-1}$ is a starting point or ending point for a vertical bridge or half
bridge.  Assume that vertex is $g_{-1}$.  (The other case is identical.)

Now, the bridge creation cannot raise or lower $g_{-1}$, since this would
necessarily change it to a vertex with edges on the right and left and, if
a bridge end, it would thus have to be the end of bridge oriented
horizontally.  It follows that $g_{-1}$ has to be either the starting
square of a bridge, the starting square of a half bridge or the ending
square of a half bridge.

For these final three cases, we will change our strategy.  Instead of
showing that these three situations can't possibly arise (in fact, they
can), we will show that whenever our algorithm does attempt to include
$g_{-1}$ in one of these roles, there is a maximal or minimal square to
the left or right of $g$ that we can use as a half bridge {\it instead} of
creating the bridge or half bridge containing $g_{-1}$.  (This is the
second rule alluded to earlier.)

First of all, in all three of these cases, the incident edges to the left
and right of $g_{-1}$ cannot be in $H$.  It follows that the edges above
the squares left and right of $g_{-1}$ must be in $\mathcal H$, as seen in
the following figure.

\begin{figure} [htb]
\begin{center}
	\includegraphics[width=3cm] {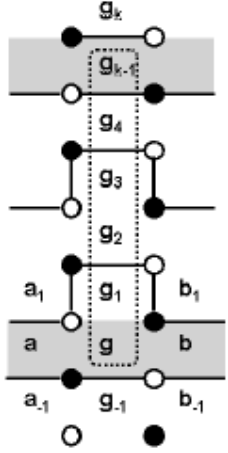}
	\caption{A bridge connects $g$ and $g_{k-1}$ when the algorithm
 attempts to include $g_{-1}$ in a bridge or half bridge.  The squares
 $a_{-1}$ and $b_{-1}$ may actually be parts of holes, but the edges
 between $a$ and $a_{-1}$ and between $b$ and $b_{-1}$ are necessarily
 present.}
\end{center}
\end{figure}

Since $g_{-1}$ was to be part of a new bridge or half bridge, it follows
that neither $a_{-1}$ nor $b_{-1}$ may be an inner bridge square.  Thus,
if there is a bridge containing $a$ (resp. $b$), it must be a bridge of
length one from $a$ to $a_{-1}$ (resp. $b$ to $b_{-1}$).  Also, we cannot
have bridges at both of these places, since they would then have to join
the same pair of boundaries, and this never happens in our algorithm.  
(If $g_{-1}$ were a triple boundary square, then $a_{-1}$ and $b_{-1}$
would not necessarily separate the same pair of cycles; however, by
assumption, this is not the case.)

Now, suppose there is a bridge from $a$ to $a_{-1}$ and not from $b$ to
$b_{-1}$.  Then we take the square $b$ to be our half bridge.  Let $B_k$
be the boundary containing $b$ --- then $M_k$ consists of the single
square $b$.  Now, we must show that $b$ cannot be adjacent to squares in
any of the other $M_i$.  We have already dealt with the case of the
squares on the left and right of $b$.

Now, suppose that $b_1$ is a bridge square on a boundary $B_l$.  This
implies that the bridge starting at $g$ must be of length two, since
otherwise, $b_1$ would be adjacent to an inner bridge point.  Since we are
joining cycles across $B_k$ with $b$, we must use the bridge containing
$g$ and $g_1$ to join cycles (at $w_1$) across $B_l$.  Thus, the bridge
starting at $b_1$ will not be used to join cycles across $B_l$, and hence
$b_1$ will not be in any $M_i$.

Similarly, suppose that $b_{-1}$ is a bridge point on a boundary $B_l$.  
Since we are joining cycles across $B_k$ with $b$, we must use the bridge
containing $a$ and $a_{-1}$ to join cycles (at $a_{-1}$) across $B_l$.  
Thus, the bridge starting at $b_{-1}$ will not be used to join cycles
across $B_l$, and hence $b_{-1}$ will not be in any $M_i$.  These figures
illustrates some of these possibilities.

\begin{figure} [htbp]
\begin{center}
	\includegraphics[width=6cm] {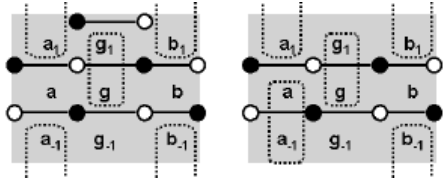}
	\caption{Examples of other bridges that may be present.}
\end{center}
\end{figure}

Now, the remaining case is when there is neither a bridge from $a$ to
$a_{-1}$ nor from $b$ to $b_{-1}$.  In this case, if there is no bridge
starting at $a_{-1}$ (resp. $b_{-1}$), we can take our half bridge vertex
to be $a$ (resp. $b$), and then the proof that this square is not adjacent
to any other $M_i$ is the same as in the previous case.

However, suppose there are bridges beginning at both $a_{-1}$ and
$b_{-1}$, and these two squares are contained on the boundary $B_l$. Then
at most one of the two can be part of $M_l$, and we can determine which
one that is from the boundary tree structure using the fact that our half
bridge will be at $B_k$.  Suppose $a_{-1}$ is a vertex that is not part of
$M_l$ (and hence not part of any $M_i$).  Then let $a$ be our final half
bridge vertex, and the proof that $a$ is not adjacent to any vertex in one
of the $M_i$ is then the same as in the previous cases.

As a final note, we have already shown that the square $f$ we use as our
final half bridge after producing a triple boundary square is not a bridge
square.  Inspection shows that the squares adjacent to $f$ across edges
contained in $H$ cannot be bridge squares since they cannot be critical
boundary squares. $\square$

\subsection{When Algorithm Gets Stuck, Number of Cycles is Minimal}

Suppose that in the above algorithm, we started with $m$ cycles, then
added $r$ new cycles and $r$ new bridges, and then got stuck.  That is, we
produced an $H$ with the following properties:

\begin{enumerate} \item $H$ divides $G$ into $m+r$ cycles, and there are
exactly $r$ non-overlapping bridges designated.

\item A square of one bridge may be adjacent to a square of another only
if the two bridges are oriented vertically (horizontally) and the two
squares are both bridge ends; no two bridges are ever adjacent at more
than one square.

\item Every critical boundary square is either a stuck square or a
square on a bridge.  No boundary squares have vertices from more than two
cycles. \end {enumerate}

Now, for each bridge $f, f_1, \ldots, f_{m-1}$, we will refer to the
squares $f_{m-2}, f_{m-1}$ as a {\it fixed pair}.  The fixed pair looks
like the length-two bridge shown in a previous diagram; that is, the two
squares are separated by an edge in $H$, the higher of the two is a local
maximum, and the lower is a local minimum.  Intuitively, these pairs occur
at places where there is slack in the height function that could be used
to join cycles together; however, we cannot make use of that slack because
we are holding the edge between those squares to be fixed.

Now, our system satisfies:

\begin{enumerate} \item $H$ divides $G$ into $m+r$ cycles, and there are
exactly $r$ designated fixed pairs which do not overlap.

\item Whenever a square of one fixed pair is adjacent to a square of
another, the two fixed pairs are both oriented vertically (horizontally).  
No two fixed pairs are ever adjacent at more than one square.

\item Every critical cycle boundary square is either a stuck square or a
member of a maximal (minimal) row of squares, whose latter endpoint is a
member of a fixed pair. \end {enumerate}

We would like to prove that any other $H'$ in the same connected component
of $S_{\phi}$ has at least $m$ distinct cycles.  We will do this by
showing that we can in a number of steps alter $H$ in such a way that each
step will preserve all the properties listed above, and so that when we
are done we will have produced a new $H'$ with potential $\tau_{H'}$, such
that $\tau_{H'}$ agrees with $\tau_{H}$ on all the squares on all but $r$
of the cycle boundaries of $H'$.  Joining the cycles across a particular
one of these boundaries can reduce the number of cycles in $H$ by at most
one, so we will conclude that $H'$ has at least $(m+r)-r = m$ separate
cycles.

Suppose that $\tau_H'$ and $\tau_H$ are not in complete agreement, and let

$$F_+ = \{f | \tau_H(f) >\tau_{H'}(f) \}$$
$$F_- = \{f | \tau_H(f) < \tau_{H'}(f) \}$$

Now, suppose that $f \in F_+$ (the $F_-$ case being symmetrical) and $f$
is maximal {\it along the boundary}.  That is, $f$ is a critical square
and $f$ is higher than its two boundary square neighbors.  Then either $f$
is maximal overall (in which case it must be a bridge square), or there is
some $f_1$ higher than $f$.  In the latter case, if $f, f_1, \ldots,
f_{m-1}$ is a maximal vertex row starting at $f$, then each of $f_1,
\ldots, f_{m-1}$ must be in $V_+$.  Since $\tau_H$ and $\tau_H'$ agree on
holes and the unbounded face, $f_1, \ldots, f_{m-1}$ is a maximal row of
squares, and hence $f_{m-1}$ is a member of a fixed pair.  Since neither
square to the left or right of $f_{m-1}$ is a local maximum or minimum,
this fixed pair must be oriented the same direction as the $f_1, \ldots,
f_{m-1}$.

\begin{figure} [htbp]
\begin{center}
	\includegraphics[width=5cm] {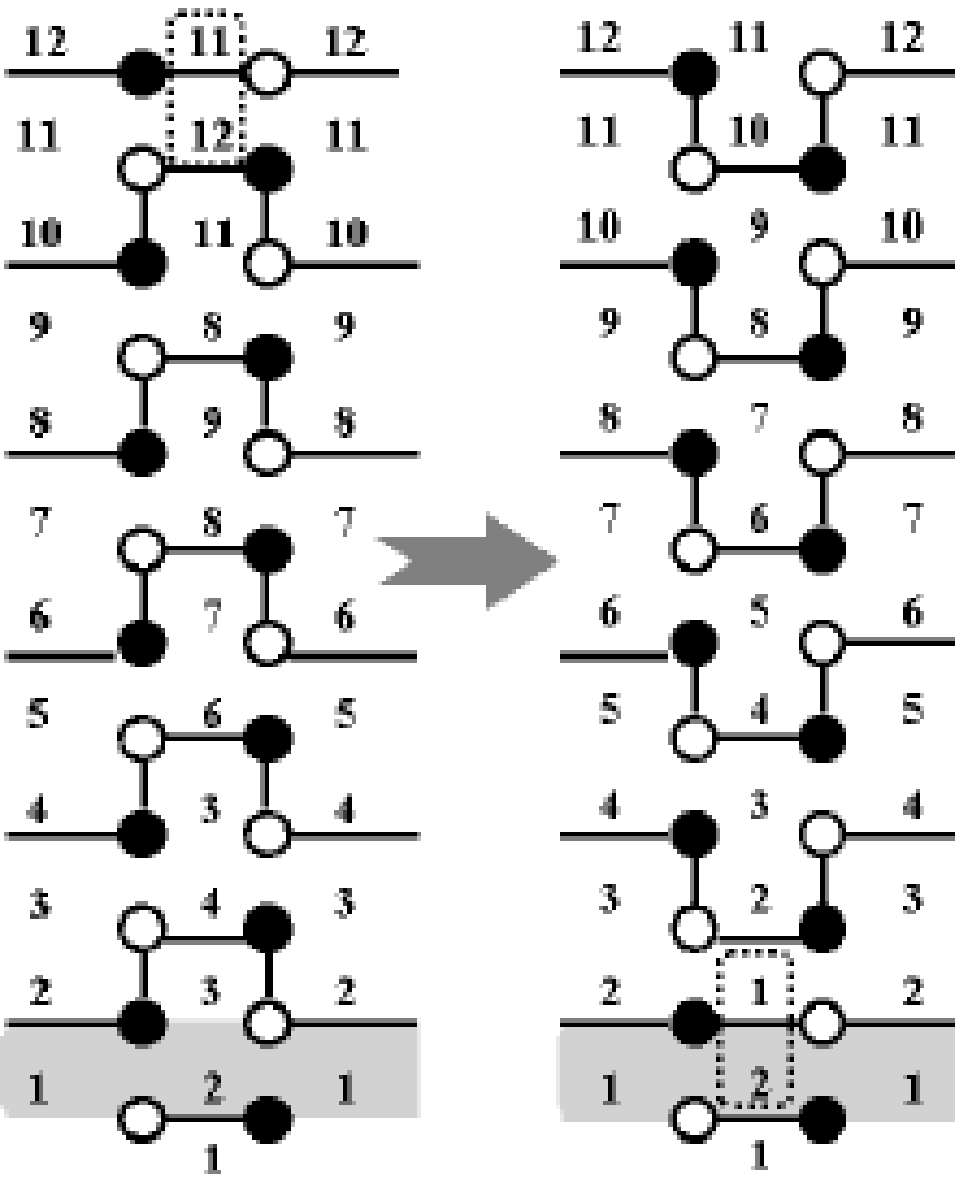}
	\caption{When a critical boundary square is part of an odd length
maximal row of length at least three and its other end square is the first
of a ``fixed pair," we can lower all the squares except the critical one
and move the fixed square from the top to the bottom.}
\end{center}
\end{figure}

Now, as the diagram shows, we will remove the fixed pair containing
$f_{m-1}$, lower the vertices $f_1, \ldots, f_{m-1}$, and add $f, f_1$ as
a new fixed pair.  One easily checks that this new system satisfies the
same three properties as the old one.  The first property is obvious.  
The second is also fairly obvious from the diagram once one notes that a
vertically (horizontally) oriented fixed pair is never incident to any
vertical (horizontal) edges.  This fact rules out a horizontally oriented
fixed pair adjacent to the new fixed pair.  A vertically oriented pair
cannot be adjacent to the new fixed pair at two places on the left or
right because one of the squares on each of the left and right is
necessarily not maximal or minimal.  The third property also follows from
inspection once we notice that each of the squares to the left or the
right of the central column in the diagram is on a directed path between
holes and/or boundary squares that were already critical (hence, newly
created critical boundary squares on the sides don't pose a problem);
also, critical squares along the row containing the new fixed pair that
were previously part of maximal or minimal rows with an endpoint in the
old fixed pair will now be part of maximal or minimal rows with an
endpoint on the new fixed pair.

We repeat this process for all other critical boundary points in $F_+$ and
similarly in $F_-$ until every critical boundary square is either

\begin{enumerate}
\item A square on which $\tau_H$ and $\tau_{H'}$ are equal.
\item A member of a fixed pair in $F_+$ that is maximal.
\item A member of a fixed pair in $F_-$ that is minimal.
\item A square in $F_+$ that is minimal along the boundary.
\item A square in $F_-$ that is maximal along the boundary.
\end{enumerate}

If $f_1$ and $f_2$ are members of a fixed pair, and $f_1$ is in $F_+$ and
is maximal, then although $f_2$ may not be in $F_+$, it cannot be in $F_-$
because it is adjacent to $f_1$.  Similarly, if $f_1$ is maximal in $F_+$,
$f_2$ cannot be minimal in $F_-$.  It follows that the number of boundary
squares in $F_+$ that are maximal along the boundary plus the number of
critical points in $F_-$ that are maximal along the boundary is at most
$r$, the fixed square pairs.

Next, let $B_i$ be a connected component of the set of boundary squares
separating some single pair of cycles.  Because there are no squares
separating three distinct cycles, one easily checks that every boundary
square must be adjacent to exactly two faces (which may or may not be
squares) on the same boundary, so that this set of squares can be written
as a sequence $f_0, \ldots, f_k$ such that each face is adjacent to its
successor and either $f_0$ and $f_k$ are adjacent to holes or they are
adjacent to each other.  In both cases, the edges between adjacent squares
are not contained in $H$.

Next, we claim that if $B_i$ has a square in $F_+$ ($F_-$, it must contain
a square in $F_+$ ($F_-$), that is maximal (minimal) along the boundary.  
If $f_i$ is the first square in $F_+$ (necessarily $f_i \succ f_{i-1}$),
then the first square after $f_i$ which is higher than its successor is a
local maximum; there must be such a place, because $\tau_H$ and
$\tau_{H'}$ agree on the endpoints.

If $B_i$ contains only squares but contains one or more squares in $F_+$
and one or more not in $F_+$, then we can apply the same argument to the
sequences of squares contained in $F_+$.  If all of the squares in $B_i$
are in $F_+$, then it is sufficient to note that every closed sequence of
boundary squares must contain both maxima and minima.  To see this, let
$f$ be any vertex on the boundary with a neighbor $f_1$ such that $f_1
\succ f$.  Then either $f_1 \prec f_2 \prec f_3 \prec \ldots \prec f$ or
there is some vertex that is a maximum along the boundary.  However, if in
fact, $f_1 \prec f_2 \prec f_3 \ldots$ then every square is a corner
square, and the sequence must progress diagonally and never return to $f$.

Thus, every connected component of the cycle boundary structure that does
not contain a vertex in $F_+$ that is maximal along the boundary or a
vertex in $F_-$ that is minimal along the boundary must have $\tau_H$ and
$\tau_{H'}$ equal everywhere.  Since there are at most $r$ vertices of
this type, there are at most $r$ connected components of the cycle
boundary structure on which $\tau_H$ and $\tau_{H'}$ disagree.  Joining
cycles across each of those $r$ components reduces the number of cycles by
at most one, so $m$ is indeed the minimum number of cycles for the
connected component of $S_{\phi}$ containing $H$.

\subsection{Runtime Analysis and Implementation}

The most naive implementation gives a runtime of at most $O(|F|^3)$ in the
worst case.  This assumes we start out with $O(|F|)$ separate cycles.  
Each time we try to join a pair of those cycles, we produce $O(|F|)$
bridges.  Each time we produce a bridge, we spend $O(|F|)$ time updating
our list of critical boundary squares.

To produce the bound $O(|F|^2)$, we still assume we start out with
$O(|F|)$ separate cycles, but we claim that given $H$, we can always
reduce $p(H)$ by one (if possible) in $O(|F|)$ steps, because it is not
really necessary to spend $O(|F|)$ steps updating our list of critical
boundary squares each time we produce a bridge.

We show this first in the polyomino case.  First, we compute $B(H)$ and
put all the critical boundary squares in a list in $O(|F|)$ steps.  Next,
each time we lower or raise a square $f$ to produce a bridge, we create a
new boundary containing $f$.  Now, if we knew that this new boundary
contained a locally maximal or minimal square, we could take that square
as our new half bridge; regardless of where it was, we could recompute
$B(H)$, find that square, and complete this phase of the algorithm in
$O(|F|)$ steps.  But for now, we will assume the boundary contains no
locally maximal or minimal squares and try to reproduce it.

That is, we attempt to form paths of non-minimal, non-maximal squares
starting at $f$ and extending in either direction such that the paths do
not pass through any edges in $H$; we stop when either the paths meet each
other (and we have a cyclic path) or they reach a previously existing
boundary square or the unbounded face.  The key observation is that under
these assumptions, there is only one way to create the path.  To see this,
suppose $f_k$ was one tentative boundary square, and $f_{k+1}$ was my next
choice, immediately north of $f_k$.  Assuming $f_{k+1}$ is not a dead end,
$H$ contains either zero, one, or two edges of $f_{k+1}$; if it contains
two edges, there is only one choice for $f_{k+2}$.  If it contains one or
zero edges, there are more choices; however, assuming $f_{k+1}$ is not
locally maximal or minimal, quick inspection shows that all but one of
these is necessarily a dead end or a local maximum or minimum.

If my sequence gets stuck in a dead end, I conclude that I must have taken
a wrong turn --- and hence my boundary contains a local maximum or
minimum, and I'm done.  If the sequence runs into a previously existing
boundary square, then either I took a wrong turn (and there's a local
maximum or minimum on the new boundary) or I've created a triple boundary
square.  If the sequence is forced to choose a minimal or maximal square,
then either that square is correct or there was another minimal or maximal
square off of a wrong turn.  In each case, there exists a half bridge
somewhere, and I can find it and complete this phase of the algorithm in
$O(|F|)$ steps.  If the sequence meets itself or if both ends reach the
unbounded face, then I conclude that I really have found the new cycle
boundary squares, and I proceed with the algorithm.

For another enhancement to this phase of algorithm, I note that each time
I create a new bridge square, it remains a bridge square until I'm
finished; I label these squares as bridge squares.  Similarly, if I find a
critical boundary square is {\it stuck}, I label the entire increasing row
of squares between that square and the other boundary square or unbounded
face edge as stuck.  If any of those become critical boundary squares at
later stages in the algorithm, I won't have to check them again.  With
this scheme, I check each square as part of a new boundary creation stage
at most once and I check each square as part of a bridge creation or stuck
critical boundary square determination step at most once.  I conclude that
this entire phase of the algorithm is $O(|F|)$.

If $G$ has holes, then when I create a new boundary, even if I don't
create any triple boundary squares, the new boundary may still contain
several connected path components, each connecting one boundary to
another.  This suggests that my previous technique of computing the
boundary squares won't work, because when my path reaches a hole, it's not
clear where to go next.

As it turns out, we can use the same idea as in the polyomino case with
minor modifications.  Suppose there are $k$ holes counting the unbounded
face.  (For the purpose of algorithm bounds, treat $k$ as a fixed
constant.)  Then we say two holes are adjacent if there is a path of
squares from one to the other that doesn't pass through any boundary
squares or cross any edges of $H$ and such that no square along the path
is maximal or minimal.  As part of our $O(|F|)$ preparation for this phase
of the algorithm, we compute $B(H)$ and determine whether each pair of
holes is adjacent, storing the path that connects them if there is one.  
(There is necessarily only one such path---if there were two, at least one
would contain a boundary face.)  We also make a table indicating when two
such paths overlap.  Observe that under this definition of adjacent, the
graph on the holes is a disjoint union of trees.

Now, after we lower or raise an $f$ to produce a bridge, we try to guess
the new boundary much the same way as before.  We produce paths of squares
extending in either direction; if they close in on each other or they both
reach the same hole, we have found our boundary.  If one reaches a dead
end or another boundary face, we conclude that there was a maximum or
minimum somewhere and we find it.  If the paths reach different holes,
then we check to see if those two paths are on the same adjacency
tree.  If they are, we take the unique set of boundaries connecting them
to be our new boundary.  If they are not, we conclude that there is a
maximum or minimum along our new boundary somewhere and we find it.

\section{Generalizing to Gridlike Graphs}

A careful examination of the steps in the above proofs reveal that the
requirement that $G$ be a grid graph is much stronger than what we needed.  
In this section, we try to capture the ``grid-like" properties that are
actually necessary.

Let $G$ be a graph embedded in on a plane or $n$-holed torus such that all
of its faces are simple polygons.  We say a face $f$ of $G$ is a {\it good
square} if it satisfies the following properties:

\begin{enumerate}
\item Every vertex of $f$ has degree at most $4$.
\item  If two adjacent vertices $v_1$ and $v_2$ of $f$ both have degree
four, then the face $f'$ which shares the edge $(v_1,v_2)$ with $f$ is a
(not necessarily good) square.
\item The square $f$ does not share more than one edge with any single
square.
\end{enumerate}

A {\it $k$-holed, $n$-holed-torus-embedded grid-like graph} is a graph $G$
embedded in an $n$-holed torus such that all of its faces are polygons and
all but $k$ of its faces are good squares.

\begin{figure} [htbp]
\begin{center}
	\includegraphics[width=5cm] {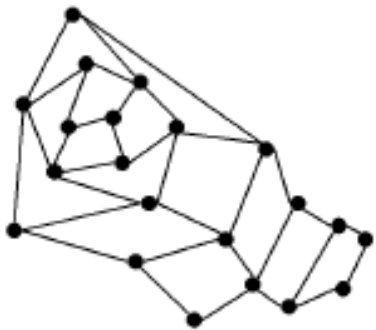}
	\caption{Planar grid-like graph with no holes.}
\end{center}
\end{figure}

In this section, we say that $H \sim H_0$ if $\omega_H - \omega_{H_0}$ is
exact and let $\tau$ be the $2$-chain with this as its boundary (such that
$\tau(f_0)=0$ for a reference hole, $f_0$).  We claim that the algorithm
developed in the previous section extends to gridlike graphs.  That is, if
we fix $n$ and $k$ to be constants, then given a fixed $H_0$ in a
$k$-holed, $n$-holed-torus-embedded grid-like graph, we can find an $H
\sim H_0$ for which $p(H)$ is minimal in $O(|F|^3)$ steps.

Both the algorithm and the proof of its effectiveness are the same as in
the grid graph case with a few modifications.  We could have introduced
this notation at the outset, but that would have made the first
read-through of the algorithm more confusing and only slightly more
enlightening.  So instead, we outline here the modifications required and
present the details as exercises.

\begin{enumerate}
\item Prove that if none of the good squares in $B(H)$ is locally maximal
or minimal, then all good squares are either corner boundary squares or
critical boundary squares.  (The first step is to use the definition of
good square to show that if two isolated vertices $v_1$ and $v_2$ on a
good square $f$ are adjacent, the face that shares the edge $(v_1,v_2)$
with $f$ is necessarily a square and hence maximal or minimal.)

\item Construct maximal or minimal rows as follows.  Let $f_0$ be a square
that is higher than its neighbors on the south, east, and west, and lower
than its neighbor to the north.  Show that if it is possible to decrease
$\tau_H(f_0)$ without changing $\tau_H$ on any holes, then $\tau_H(f_0)$
must be part of a unique increasing path of good squares such that every
square in the path is higher than every face outside the path.  For $i>0$,
define the south edge of $f_i$ to be the single edge that borders
$f_{i-1}$, and define north, east, and west edges accordingly.  Show that
each $f_i$ in the sequence is higher than its neighbors to the east and
west.

\item In what ways can squares in a maximal or minimal row be adjacent to
other squares in the same row?  Answer this question with the following
steps.  Let $i$ be the first integer such that $f_i$ shares an edge $e$
with some earlier $f_j$ in the sequence (not counting its neighbor to the
south).  Note that $e$ cannot be the east or west edge of $f_i$ because if
it were, $f_i$ would have to be higher than $f_j$, and $f_j$ is higher
than all squares except for its neighbor to the north.  Suppose $e$ were
the north edge of $f_i$ and the west (or, symmetrically, east) edge of
$f_j$.  If $f'$ is the neighbor of $f_j$ to the west, then $f', f_j, f_i,
f_{i+1}$ is an increasing sequence of squares sharing the same vertex, and
this also cannot happen.  Hence, it must be the case that the north edge
of $f_i$ is the south edge of $f_j$, and thus $j=0$, since already we know
the southern neighboring face of the other $f_j$.

In this case, we refer to the sequence $f_0, \ldots, f_j$ as a {\it cyclic
maximal good square row}.  Since $G$ is bipartite, it has even length
(i.e., $j$ is odd).  It is impossible to change the value of $\tau_H$ on
these good squares with $Z$-transformations alone, but we can decrease the
value of $\tau$ on all the squares in this row at once.  Show that this
has no effect on the value of $p(H)$.

\item Show that if $H \sim H_0$, it is possible to move from $H$ to $H_0$
by successively applying $Z$-transformations and raising (lowering) cyclic
minimal (maximal) good square rows.  As before, this is done by finding an
element in $F_+$ (or $F_-$) that is higher (lower) three squares.  We saw
earlier that this is possible provided $F_+$ borders at least one element
not in $F_+$.  Show that this fails only when $F_+$ contains every face in
$F$ and there are no holes and every face is higher than exactly two
squares.  (Certain cycle partitions on $C_m \times C_n$ have this
property.  Show that in this exceptional case, $S_{\phi}$ contains exactly
two elements.)

\item Show that all of the steps in the basic algorithm are well-defined
in this new setting.  Adapt the proof of the lemma that says the maximal
and minimal rows used to join vertices are disjoint and no element of one
is adjacent to any element of another.  One step requires showing that if
$f_0, \ldots, f_{m-1}$ is a north-south oriented bridge, the faces
adjacent to $f_1, \ldots, f_{m-2}$ on the left and right are on directed
paths from one hole or member of $B_H$ to another; in our new setting, we
have to deal with the case that a single square $f$ is adjacent to two or
more of the $f_1, \ldots, f_{m-1}$ on the east or west.  Still, we can let
$f_0', f_1', f_2`, \ldots, f_{m-1}, f_m'$ be the faces adjacent to $f_0,
\ldots, f_m$ and show that even if the list contains repeats, it still
describes an increasing path from one hole or boundary face to another.

\item Using fixed pairs, show that when the algorithm fails, $p(H)$ is
minimal.

\end{enumerate}

\section{Markov Chains for Hamiltonian Cycles}

Write $S_{\phi}^k = \{ H\in S_{\phi} | p(H) = k\}$, i.e., the set of
partitions of $G$ into exactly $k$ cycles.  In this section, instead of
simply using $Z$-transformations (raising or lowering single squares), we
say two subgraphs of $S_{\phi}$ are adjacent if one can be obtained from
the other by raising (lowering) a minimal (maximal) row.  Let $\mathcal H$
be a chosen homology class of $S_{\phi}$.

\begin{theorem}
Under the above definition of adjacency, $(S_{\phi}^1 \cup S_{\phi}^2)
\cap \mathcal H$ is connected whenever both $S_{\phi}^1 \cap H$ and
$S_{\phi}^2 \cap H$ are nonempty.
\end{theorem}

To see this, we show that given an $H' \in S_{\phi}^1$, and an $H \in
S_{\phi}^1 \cup S_{\phi}^2$ with $H \sim H'$, we can always successively
raise or lower rows of $H$ in a way that decreases $|\tau_{H} -
\tau_{H'}|$.

The first observation is that if $H \in S_{\phi}^1$, then we can choose
any square $f$ in $G$ on which $\tau_{H}$ is too high (too low) such that
$f$ is higher (lower) than three or four of its neighbors.  Then we lower
(raise) the maximal (minimal) row beginning with $f$, thereby decreasing
$|\tau_{H} - \tau_{H'}|$ by the length of the row.

Second, suppose $H \in S_{\phi}^2$.  Then because $H'$ does not contain
any boundary vertices, there exists at least one element of $B(H)$ on
which $\tau_H$ and $\tau_{H'}$ do not agree; assume without loss of
generality that there is an $f\in B(H)$ with $\tau_H(f)>\tau_{H'}(f)$.  
Then there must be at least one square $f$ in $B(H)$ that is maximal along
the boundary.  If $f$ is higher than three or four of its neighbors, we
lower the maximal row beginning with $f$; if the length of the row is odd,
the result is in $S_{\phi}^1$.  If the length is even, the result remains
in $S_{\phi}^2$.

We have to be a little more careful in the event that $f$ is higher than
two neighbors which are boundary faces and lower than the other two,
as in the following diagram:

\begin{figure} [htbp]
\begin{center}
	\includegraphics[width=3cm] {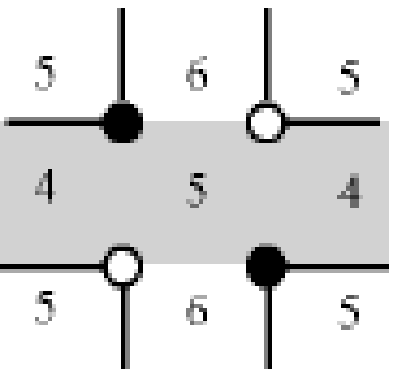}
	\caption{}
\end{center}
\end{figure}

In this case, there is a maximal row beginning with $f_1$ (the face above
$f$ in the diagram) extending upward and another beginning with $f_{-1}$
extending downward, both of which must be lowered before we can decrease
the value of $\tau_H$ on $f$ itself.

Now, if $f_1$ or $f_{-1}$ is the beginning of an even length maximal row,
we can pull that down without changing the number of cycles in $H$.  If
either square is the beginning of an odd length maximal row of length
greater than one, we can first pull up one of the squares adjacent to $f$
on the left or right in the diagram (which joins cycles and moves $H$ into
$S_{\phi}^1$) and subsequently lower the maximal row beginning with $f_1$
or $f_{-1}$ (which moves $H$ back into $S_{\phi}^2$).  Raising the square
on the left or right of $f$ is a step in the wrong direction---increasing
$|\tau_H - \tau_{H'}|$ by one---but it allows us to take three or more in
the right direction.

Suppose, however, that both $f_1$ and $f_{-1}$ are local maxima.  Then
raising a square to the left or right of $f$ and subsequently lowering one
of $f_1$ or $f_{-1}$ does not change the value of $|\tau_H - \tau_{H'}|$.  
However, suppose the longer of the two cycles contains $j$ vertices and
the other contains $k$; then we can first raise one of the squares to the
left or right of $f$ and then lower whichever of $f_1$ or $f_{-1}$ belongs
to the smaller of the cycles; this has the effect of decreasing the length
of the smaller cycle.  Since the length of the smaller cycle cannot
decrease indefinitely, if we repeat this process, we must eventually reach
a position at which it is possible to decrease the value of $|\tau_H -
\tau_{H'}|$. $\square$

Similar arguments can be made for the connectedness of $(S_{\phi}^k \cup
S_{\phi}^{k+1})\cap \mathcal H$ under these operations.

Now that we know that $(S_{\phi}^1 \cup S_{\phi}^2) \cap \mathcal H$ is a
connected graph, we can take a random walk on this graph.  To do this, we
define a transition step as follows:

Choose a face $f$ at random.  If $H \in S_{\phi}^1$, then if $f$ is higher
(lower) than three or four neighboring squares and it is possible to lower
(raise) the maximal (minimal) row beginning with $f$, then do so.  
Otherwise do nothing.

If $H \in S_{\phi}^2$, then if $f$ is in $B(H)$ and $f$ is higher (lower)
than three or four neighboring squares and it is possible to lower (raise)
the maximal (minimal) row beginning with $f$, then do so.  Otherwise, do
nothing.  If $f \not \in B(H)$, but $f$ is the beginning of an even length
maximal (minimal) row such that lowering (raising) that row leaves $H$ in
$S_{\phi}^2$, then do so.  Otherwise, do nothing.

Because this transition process describes a simple random walk on a graph
in which each $H$ has (counting self-loops) degree $|F|$, it has uniform
stationary distribution.  Hence, applying this transition process enough
times allows us to sample randomly from the space $S_{\phi}^1 \cup
S_{\phi}^2$.  The drawback is that we do not yet have rigorous bounds on
the rate of convergence of this Markov Chain.  If $S_{\phi}^1\cap \mathcal
H$ is a reasonable fraction of $(S_{\phi}^1 \cup S_{\phi}^2) \mathcal H$,
we can perform repeated Markov searches until we find an element of
$S_{\phi}^1$; this gives us a way of choosing a Hamiltonian cycle of $G$
at random.  A second drawback is that it is not known how ``reasonable" a
fraction $S_{\phi}^1 \cap \mathcal H$ occupies of the entire space
$(S_{\phi}^1 \cup S_{\phi}^2) \cap \mathcal H$.

Nonetheless, we suspect that the fraction and the convergence rate are
inverse polynomial and polynomial respectively in the number of faces of
$G$; we have shown this to be the case for some special families of grid
graphs.

\section {Almost Taut Examples}

A chipped Aztec diamond of size $2n$ is created as follows.  Start with an
Aztec Diamond graph of radius $2n$.  Then successively remove ``chips"
from the upper half of the diamond; here, a {\it chip} is a horizontally
adjacent pair of vertices $v_1, v_2$ with a white vertex on the left such
that both vertices are incident to the unbounded face and none of the
vertices above $v_1$ and $v_2$ remains in the graph.  One can prove that
an Aztec diamond of size $2n$ with $k<n$ chips can be partitioned into at
fewest $n-k$ disjoint cycles.  (The proof follows the lines used in the
proof of the correctness of our algorithm: that is, we show inductively
that we can separate the vertices into $n$ distinct cycles with $k$ fixed
pairs in a configuration satisfying the three properties used in our
proof.)  Hence $k$ must be at least $n-1$ in order for the diamond to have
a Hamiltonian cycle.  We will assume in the remainder of this discussion
that $k$ is exactly $n-1$.  The following is an example of a Hamiltonian
cycle on a chipped Aztec diamond with $n=4$.

\begin{figure} [htbp]
\begin{center}
	\includegraphics[width=7cm] {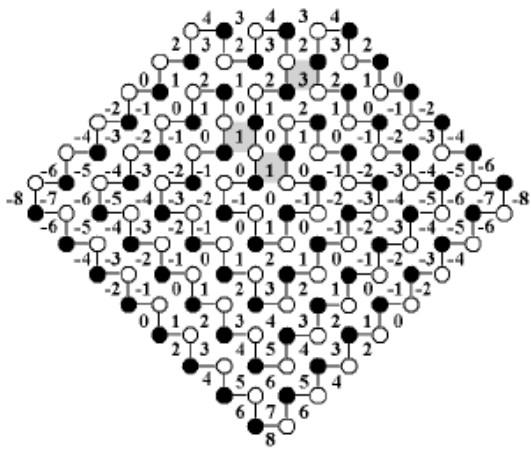}
	\caption{Chipped Aztec Diamond}
\end{center}
\end{figure}

The values of $\tau_H$ on the lower half of the diamond are all fixed.  
Now, label the squares with coordinates so that the center squares is
$(0,0)$ and consider the $4n-1$ vertical columns of squares in the upper
half, each beginning with a square of $y$-coordinate zero.  In a sequence
of squares $f_0, \ldots, f_i$ from the bottom of such a column to the top,
the number $k$ for which $f_k \prec f_{k-1}$ is determined by the boundary
data.  Let $a_H(i,j)$ be the value $k$ for which such an event occurs in
the $i$th column for the $j$th time.  The set of $i$ and $j$ for which
$a_H(i,j)$ is defined is independent of $H$.  In the example, we have
colored in the three squares in the upper half on which $\tau$ is higher
than it is on the square immediately to the north.  Note that the middle
square cannot take a value higher than than the outer two.  That is,
$a_H(0,1)$ is less than or equal to $a_H(1,1)$ and $a_H(1,0)$.  We can say
that the set of possible $H$ corresponds to the set of labelings of an
L-shaped Young tableaux with three squares with numbers between one and
five such that the numbers are non-decreasing in each row or column.  The
reader may deduce a more general correspondence between chipped Aztec
diamond cycle partitions and Young tableaux labelings.  There is also a
more subtle correspondence between Hamiltonian cycles of a chipped Aztec
diamond and {\it standard} labelings of Young tableaux.

This next example is even simpler than the last.

\begin{figure} [htbp]
\begin{center}
	\includegraphics[width=7cm] {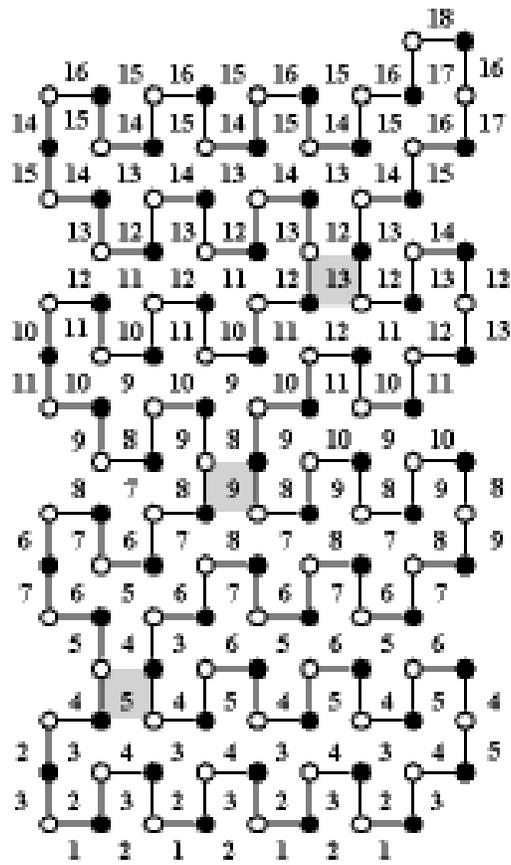}
	\caption{Low Slack Tower}
\end{center}
\end{figure}

Here there are only three squares on which $\tau_H$ takes a higher
value
than it does on the square immediately to the north, and each necessarily
occupies on of the three columns the squares are in.  The reader may show
that the set of Hamiltonian cycles corresponds to the permutations of
three elements and deduce a more general result.  It is a worthwhile
exercise to work out the relative sizes of $S_{\phi}^1$ and $S_{\phi}^2$
and try to put upper bounds on the Markov chain convergence rates for
these simple examples.  Intuitively, one would expect that in examples
that are not so almost-taut, a randomly chosen element of $S_{\phi}$ would
more likely to have one and perhaps several boundary squares beginning
maximal or minimal rows of odd length (unlike the above two examples,
in which most would have even length).  Hence, one might expect the
difference in sizes between $S_{\phi}^1$ and $S_{\phi}^2$ to be less
pronounced in these cases.  It is not yet known how to make this rigorous.

\section{Conclusion and Open Questions}

The use of height functions for Hamiltonicity computation and sampling
problems is a new field of research; there are doubtlessly many
improvements to be made, empirical studies to be done, bounds to be
determined, and generalizations to be drawn.  There also remain many open
problems regarding mean properties of restricted difference height
functions and performance of Coupling from the Past algorithms.  We list
just a few possible directions for further study.

\begin{enumerate}

\item Can we bound the Markov chain convergence for the set of Hamiltonian
cycles of a grid graph?  Can we generalize the Hamiltonian cycle counting
algorithms and bounds given in \cite{BT}, \cite{KZ}, \cite{KR}, and
\cite{SS} for rectangular grids?

\item Can the log of the number of Hamiltonian cycles in a large simply
connected grid graph be approximated with an entropy integral the way the
number of restricted vertex degree subgraphs can be (see \cite{CKP})?

\item A close look at the grid graph algorithm shows that we have actually
produced a characterization of grid graphs that have cycle covers and fail
to be Hamiltonian; they are those for which ``fixing" some $r$ edges in
$H$ pulls the height function so taut that all the cycle boundary squares
separating $H$ into a set of $m+r$ regions are held in place (for some
$m>0$).  This suggests there may be a faster algorithm for recognizing,
for example, polyomino graphs with this property.  Is there a linear or
$O(|F|log|F|)$ algorithm for recognizing Hamiltonian polyomino graphs?

\item If we assign a weight to each edge in $G$, there are well-known
constrained network optimization algorithms that will enable us to find an
$H$ partitioning $G$ into cycles such that the sum of the edges in $H$ has
minimal weight.  What can we say about the problem of finding a
Hamiltonian cycle of minimal weight?  Will the techniques of this paper
help in this special case of the traveling salesman problem?

\item Is there a fast algorithm for computing the {\it maximal} degree of
elements of $S_{\phi}$?  (In some sense, this maximal degree quantifies
the notion of ``slack" in the height function.)

\item Can we use techniques similar to those for grid graphs to determine
the Hamiltonicity of hexagonal graphs?  (It seems at least plausible;
there are nice analogs to maximal and minimal square rows.  The tricky
issue is parity---a $Z$-transformation always changes the number of cycles
by zero or two.  Thus, to decrease the number of cycles, one must find a
$Z$-transformation that will join three cycles together into one.)  If I
can compute the Hamiltonicity of hexagonal graphs with these techniques,
can I go still further?  Is there a larger class of graphs for which these
or similar techniques give polynomial algorithms for determining
Hamiltonicity or other graph invariants?

\end{enumerate}

\section{Acknowledgements} Thanks to James W. Cannon for suggesting the
problem and to Joe Gallian and others at Duluth.  Thanks also to Persi
Diaconis, Amir Dembo, Robin Pemantle, and Alan Hammond for helpful
conversations.  Thanks to my wife, Julie, for providing not only love and
support but several nice diagrams as well.

\end{document}